\documentclass[11pt,oneside, final]{amsart}
\copyrightinfo{0}{Iranian Mathematical Society}
\pagespan{1}{\pageref*{LastPage}}
\usepackage{etoolbox,lastpage}
\commby{}
\usepackage{amsmath,amsthm,amscd,amsfonts,amssymb,enumerate}
\usepackage{graphicx}
\usepackage{color}
\usepackage[colorlinks]{hyperref}

\makeatletter \@addtoreset{equation}{section}

\newtheorem{theorem}{Theorem}[section]
\newtheorem{proposition}[theorem]{Proposition}
\newtheorem{lemma}[theorem]{Lemma}

\theoremstyle{definition}

\newtheorem{remark}[theorem]{Remark}

\theoremstyle{remark}
\numberwithin{equation}{section}
\begin{document}

%% The title of the paper goes here.  Edit your title.

\title[Legendrian curve flow in Sasakian sub-Riemannian 3-manifolds]{Legendrian curve flow in Sasakian sub-Riemannian 3-manifolds$^{**}$}
%% Now edit the following to give First Author name and address:
%% $^*$ for the corresponding author.

\author[J. Cui]{Jingshi Cui}
%%\address[Jingshi Cui]{School of Mathematical Sciences, Xiamen University, Xiamen 361005, P. R.China}\email{cuijingshi626@163.com}

\author[P. Zhao]{Peibiao Zhao$^*$}
%%\address[Peibiao Zhao]{Dept. of Applied Mathematics, Nanjing University of Science and Technology, Nanjing 210094, P. R.China} \email{pbzhao@njust.edu.cn}
%% If there are three of more authors they are added in the obvious way.

\thanks{$^*$Corresponding author}
%%\thanks{$^{**}$ Supported by NNSF of China(No.12271254; No.12141104)}
%------------------------------------------------------------------------------------%
%%
%% Use the following command to make the title for the paper.
%
 %\CoverPage
\begin{abstract}
    In this paper, we introduce a kind of inverse mean curvature flow (\ref{1.2}) in a Sasakian sub-Riemannian 3-manifold $M$ for Legendrian curves, which slightly differs from the classical one, and confirm that this flow preserves the Legendrian condition and increases the length of curves. We establish the long-time existence of the flow (\ref{1.2}) when the Webster scalar curvature $W$ of $M$ satisfies $ W \in (-\infty,  \bar{W}_{0} )\cup \{ 0\} \cup (W_{0}, +\infty)$, where $\bar{W}_{0} <0$ and $W_{0} >0$ are constants. Moreover, we derive that the local limit curve (the asymptotic behavior) along the flow (\ref{1.2}) is a geodesic of vanishing curvature when $W \geq 0$, wherea it is a geodesic of nonvanishing curvature when $W$ is a negative constant.

    Specially, in the first Heisenberg group $\mathbb{M}(0)$, we further construct a length-preserving flow (\ref{1.3}) via a dilation of the flow (\ref{1.2}) and show that closed Legendrian curves converge to Euclidean helices with vertical axis. By exploiting the properties of the flow (\ref{1.3}), we establish a Minkowski-type formula for Legendrian curves in $\mathbb{M}(0)$ and provide a new proof of the fact that the total curvature of $\gamma \subset \mathbb{M}(0)$ with strictly positive curvature equals $2\pi$.
\end{abstract}

\maketitle
%
%%% The following environment is needed for the abstract.
%%%

\section{Introduction}
Let $(M,g)$ be a Riemannian manifold, and $X:\Sigma \times[0,T_{0}) \to M$ be a family of smoooth immersions $X_{t}(\cdot) = X(\cdot, t)$ of $\Sigma $ into $M$. The classical form of the inverse mean curvature flow (IMCF) is given by 
\begin{align}\label{1.1}
    \frac{\partial X}{\partial t}(p,t) =\frac{\nu }{H}
\end{align}
where $H$ and $\nu$ denote the mean curvature and the outward unit normal of $\Sigma_{t}=X_{t}(\Sigma)$. The study of the flow (\ref{1.1}) has encompassed both classical solutions and weak solutions. For compact, star-shaped hypersurfaces with strictly positive mean curvatures, it was proved by Gerhardt \cite{G1990} and Urbas \cite{U1990} that the smooth solutions evolving under IMCF exist for all time and converge to round spheres after a suitable rescaling. The topological sphere hypersurfaces in Euclidean space to expand under the IMCF is a prototypical behavior of the flow. However, on general Riemannian manifolds or starting from general initial hypersurfaces, as stated in \cite{H2019}, the smooth IMCF inevitably forms a finite-time singularity. To address this situation, Huisken and Ilmanen \cite{HI2001} introduced a notion of weak IMCF by finding a suitable variational principle for the level set equation of (\ref{1.1}).

The inverse mean curvature flow, in both its smooth and weak formulations, has found numerous applications. For instance, the smooth IMCF is frequently used to establish geometric inequalities, as demonstrated in \cite{BHW2016,GR2020,KWWW2022,KW2023}. The weak IMCF has been employed in establishing the Riemannian Penrose inequality \cite{HI2001}, computing Yamabe invariants \cite{BN2004}, investigating problems related to scalar curvature \cite{BM2008,LN2015,S2016}, and providing a new proof of the Ricci pinching conjecture \cite{HK2024}. These applications are primarily based on finding monotone quantities along the flow and analyzing the asymptotic behavior of the (rescaled) limit. Recent studies on the properties and applications of the IMCF have achieved significant progress, and the research has gradually extended to sub-Riemannian spaces.
 
Sub-Riemannian spaces are degenerate Riemannian manifolds in which the inner product is defined only on a subbundle of the tangent bundle, with the Heisenberg group providing a simplest nontrivial model. In this setting, the natural analogue of the classical IMCF is the horizontal inverse mean curvature flow (HIMCF). Due to the presence of characteristic points in the Heisenberg group, where the tangent hyperplane coincides with the horizontal distribution, the level-set formulation of the HIMCF becomes degenerate, and consequently only weak solutions can be considered in general. Inspired by the ideas in \cite{HI2001,M2007}, Cui and Zhao \cite{CZ2023} introduced weak solutions of the HIMCF as minimizers of energy functionals and established their existence and uniqueness in the first Heisenberg group. The existence of a global weak solution of the HIMCF in the $(2n+1)$-dimensional Heisenberg group was later established by Pisante and Vecchi \cite{PV2024}. Nevertheless, analyzing the asymptotic behavior of the HIMCF remains a substantial challenge. The main difficulty arises from the lack of a sufficiently developed differential geometric framework for hypersurfaces in sub-Riemannian spaces. In addition, the nonzero dimensions and the complex structure of the set of characteristic points create further difficulties that do not appear in the Euclidean setting.

This paper investigates the inverse mean curvature flow on Sasakian sub-Riemannian 3-manifolds, which are pseudohermitian manifolds with vanishing torsion. According to classical results \cite{FG1996,T1969}, the sub-Riemannian 3-space form $\mathbb{M}(K)$ is the unique complete, simply connected Sasakian sub-Riemannian 3-manifold with constant Webster scalar curvature $K/2$. In particular, the first Heisenberg group corresponds to the case $K = 0$, which we denote by $\mathbb{M}(0)$ for convenience. It is worth emphasizing that Sasakian geometry is closely related to various areas of mathematics and physics, and a comprehensive introduction with numerous examples can be found in \cite{CK2008}

In Sasakian manifolds, there exists a class of curvature flows that preserve the Legendrian condition and exhibit favorable properties concerning long-time existence, convergence, and applications. Let $\Sigma^{n}$ be a smooth $n$-dimensional manifold and $F: \Sigma^{n} \times [0, T_{0}) \to M$ be a family of immersions into a Sasakian $(2n+1)$-manifold $(M,\theta,g_{\theta},\nabla, J)$. For each $t$, denote $F_{t}:= F(p,t)$ for $p \in \Sigma^{n}$ and $(F_{t})_{k}:= \frac{\partial F_{t}}{\partial p_{k}} \in T_{F_{t}(p)} M$. Smoczyk \cite{S2003} first introduced the following evolution equation
\begin{align*}
    \frac{\partial }{\partial t} F =  \nabla^{k} f \nu_{k}  +  2f T
\end{align*}
where $\nu_{k}$ are defined by $\nu_{k} = J (F_{t})_{k} $, $f : \Sigma^{n} \times [0, t_{0}) \to \mathbb{R}$ is a smooth family of functions and $T$ is the Reeb vector field. A typical example of the above evolution is the Legendrian mean curvature flow (LMCF), which not only decreases the volume energy as in the classical mean curvature flow but also preserves the Legendrian condition along the flow. Smoczyk \cite{S2003} investigated the LMCF on pseudo-Einstein Sasakian manifolds, and in the one-dimensional case the flow reduces to the Legendrian curve shortening flow, which has been applied to the isotopy problem for Legendre curves. Drugan, He, and Warren \cite{DHW2018} studied this flow in $\mathbb{R}^{3}$ with the standard contact structure and obtained a partial answer to Grayson’s conjecture \cite{Gr1989}. Subsequently, Pan and Sun \cite{PS2021} proposed a different curve shortening flow in Sasakian 3-manifolds, which coincides with Smoczyk’s flow only in the first Heisenberg group, and established results on long-time existence, convergence, and the formation of Type I singularities. More recently, Chang, Han, and Wu \cite{CHW2024} proved the existence of the long-time solution and asymptotic convergnce along the LMCF in the $\eta$-Einstein Sasakian $(2n+1)$-manifolds. These contributions demonstrate the significance of Legendrian curvature flows within Sasakian manifolds.

In view of the significant progress on the Legendrian curvature flow and its favorable geometric properties, this paper investigates the curve case of the Legendrian inverse mean curvature flow in Sasakian sub-Riemannian 3-manifolds. Let $\gamma: S^{1} \times [0,T_{0}) \to M$ be a closed curve immersed in $M$ and $s$ be the arc-length parameter. We consider the following flow: 
\begin{equation}\label{1.2}
    \begin{cases}
        \frac{\partial \gamma}{\partial t} = -k^{-1} J \dot{\gamma} - 2\int_{0}^{s} k^{-1} ds T  ,\\
    \gamma(\cdot ,0) = \gamma_{0}.
    \end{cases}
\end{equation}
where $k$ is the curvature of $\gamma$ and $\dot{\gamma} = \frac{\partial \gamma}{\partial s}$ is the unit tangent vector of $\gamma$. The flow (\ref{1.2}) differs from the one-dimensional case of the IMCF by the addition of a second term on the right-hand side. This modification is designed to preserve the Legendrian condition and follows the approach in \cite{S2003}. We refer to the flow (\ref{1.2}) as the Legendrian curve expanding flow, since under this flow the length of curves increase exponentially (see Proposition \ref{prop3.1}), resembling the behavior of IMCF in Euclidean space. Firstly, we prove the long-time existence of the flow (\ref{1.2}) by establishing uniform bounds on the higher order derivatives of the curvature.

\begin{theorem}\label{thm1.1}
    Let $\gamma_{0}$ be a closed Legendrian curve in Sasakian sub-Riemannian 3-manifold $M$, and let $\gamma (\cdot, t)$ be a solution to the Legendrian curve expanding flow  (\ref{1.2}) on the time interval $[0, T_{0})$. Assume that the Webster scalar curvature $W$ of $M$ satisfies
    \begin{align*}
        W \in (-\infty,  \bar{W}_{0} )\cup \{ 0\} \cup (W_{0}, +\infty),
    \end{align*} where $\bar{W}_{0} <0$ and $W_{0} >0$ are constants. Then, there exists $\varepsilon > 0$ such that $\gamma (\cdot, t)$ exists and is smooth on the extended time interval $[0, T_{0}+\varepsilon)$.
\end{theorem}
\begin{remark}\label{rem1.2}
   In \cite{PS2021}, the long-time existence and convergence of the curve shortening flow in Sasakian 3-manifolds were established under the assumptions that the Webster scalar curvature $W$ of $M$ is bounded above by a negative constant $C_{0}$ and that the curvature $k$ of the initial curve $\gamma_{0}$ satisfies $k^{2} + 2C_{0} < 0$. In this paper, we construct different auxiliary functions in the cases $W\geq W_{0} >0$, $W\leq \bar{W}_{0}<0$ and $W=0$ to obtain higher order curvature estimates, thus establishing the long-time existence of the flow (\ref{1.2}).  
\end{remark}

For Sasakian 3-manifolds with nonnegative Webster scalar curvature, it is possible to analyze and characterize the asymptotic behavior of the flow (\ref{1.2}).
\begin{theorem}\label{thm1.3}
    Suppose that the Webster scalar curvature $W$ of $M$ satisfies $W\geq 0$ and $\gamma (\cdot, t)$ is a solution to the Legendrian curve expanding flow  (\ref{1.2}) on the time interval $[0, \infty)$. Then, for every $l \geq 0$, $\lim_{t \to \infty} | \partial_{s}^{l} k|^{2} = 0$. Moreover, if $M$ is complete, for any compact set $\mathcal{K}  \subset M$, the local limit curve $\gamma_{\mathcal{K}}(\cdot,\infty) : = \lim_{t\to \infty} (\gamma(\cdot, t) \cap \mathcal{K})$ exists and is a geodesic of vanishing curvature.
\end{theorem}

In particular, the flow (\ref{1.2}) exhibits distinct asymptotic behavior in Sasakian 3-manifolds with negative constant Webster scalar curvature.
\begin{theorem}\label{thm1.4}
    Suppose that the Webster scalar curvature $W$ of $M$ is a negative constant and $\gamma (\cdot, t)$ is a solution to the Legendrian curve expanding flow  (\ref{1.2}) on the time interval $[0, \infty)$. Then, $\lim_{t \to \infty} k^{2} = -2W$. Moreover, the local limit curve $\gamma_{\mathcal{K}}(\cdot,\infty)$ exists and is a geodesic of nonvanishing curvature.
\end{theorem}
\begin{remark}\label{rem1.5}
    Based on the results concerning geodesics in the model spaces of the sub-Riemannian 3-space forms $\mathbb{M}(K)$ with $K =-1,0,1$ (see Lemmas \ref{lem2.5} and \ref{lem2.6}), we obtain the following conclusions:  
    (i) In $\mathbb{M}(1)$, $\gamma_{\mathcal{K}}(\cdot,\infty)$ is diffeomorphic to a circle;  
    (ii) In $\mathbb{M}(0)$, $\gamma_{\mathcal{K}}(\cdot,\infty)$ is a horizontal straight line;  
    (iii) In $\mathbb{M}(-1)$, the vector field $V_{\infty}$ associated with the variation of $\gamma_{\mathcal{K}}(\cdot,\infty)$ satisfies $g_{\theta}(V_{\infty},T)(s)=\frac{1}{2}s^{2}$.
\end{remark}

In the rest of this paper, we study the Legendrian curve expanding flow (\ref{1.2}) in the first Heisenberg group $\mathbb{M}(0)$. When $\gamma \subset \mathbb{M}(0)$ evolves under the flow (\ref{1.2}), its projection onto the $xy$-plane evolves according to the IMCF of planar curves, and their curvatures coincide. Furthermore, by applying the Heisenberg dilation with factor $\lambda = e^{-t}$ to the solution $\gamma$ of (\ref{1.2}), we obtain a family of rescaled curves $\tilde{\gamma}$ that satisfy the following equation:
\begin{equation}\label{1.3}
    \begin{cases}
        \frac{\partial \tilde{\gamma}}{\partial t} = -\left(k^{-1} + g_{\theta}(\tilde{\gamma}, J\dot{\tilde{\gamma}})\right) J\dot{\tilde{\gamma}} - 2\left(\int_{0}^{s} k^{-1} ds  + g_{\theta}(\tilde{\gamma}, T) \right)T,  \\
        \tilde{\gamma}(\cdot ,0) = \tilde{\gamma}_{0}.
    \end{cases}
\end{equation}

For the rescaled flow (\ref{1.3}), we prove the following results.
\begin{theorem}\label{thm1.6}
    Let $\tilde{\gamma}_{0}$ be a smooth closed Legendrian curve in $\mathbb{M}(0)$. Then the initial value problem (\ref{1.3}) admits a smooth solution $\tilde{\gamma}(\cdot,t)$ for all $t \in [0,\infty)$. Moreover, the flow (\ref{1.3}) preserves the Legendrian condition, and the length of $\tilde{\gamma}(\cdot,t)$ remains to be a constant. As $t \to \infty$, $\tilde{\gamma}(\cdot,t)$ convergence to a Euclidean helix with vertical axis. 
\end{theorem}

Obviously, the rescaled flow (\ref{1.3}) is a kind of length-preserving flow. By exploiting its geometric properties, we establish a Minkowski-type formula for Legendrian curves in $\mathbb{M}(0)$, which provides a useful integral formula for studying curvature flows in $\mathbb{M}(0)$. 
\begin{proposition}\label{prop1.7}
    Assume  that $\gamma$ is a smooth closed Legendrian curve in $\mathbb{M}(0)$. Then,
    \begin{align*}
        \int_{\gamma} 1+ k g_{\theta}(\gamma, J \dot{\gamma}) ds = 0.
    \end{align*}
\end{proposition}

Moreover, we obtain a global property for Legendrian curves in $\mathbb{M}(0)$. Let $\gamma$ be a simple, closed Legendrian curve, which can be parameterized by arc-length $s$. The total curvature of $\gamma$ is defined as 
\begin{align*}
    \int_{\gamma} |k(s)| ds.
\end{align*} 
Using the convergence results of the flow (\ref{1.3}), we obtain the following result.
\begin{proposition}\label{prop1.8}
    Let $\gamma$ be a smooth closed Legendrian curve in $\mathbb{M}(0)$ with strictly positive curvature. Then,
    \begin{align*}
        \int_{\gamma} k ds = 2\pi.
    \end{align*}
\end{proposition}
\begin{remark}\label{rem1.9}
    In the Heisenberg group, some literature refers to Legendrian curves as horizontally regular curves, for example, in \cite{CHL2017,CFH2018,CH2019}. The Fenchel-type theorem of Chiu and Ho \cite{CH2019} states that for any closed horizontally regular curve,
    \begin{align}\label{1.4}
        \int_{\gamma} k ds \geq  2\pi.
    \end{align}
    Denote $\gamma^{*}$ be the the projection of $\gamma$ onto the $xy$-plane. Equality in (\ref{1.4}) holds if and only if $\gamma^{*}$ is a simple convex curve. By Lemma \ref{lem7.1}, $\gamma$ and $\gamma^{*}$ have the same curvature, so the convexity of $\gamma^{*}$ is equivalent to $\gamma$ having strictly positive curvature. Therefore, Proposition \ref{prop1.8} corresponds to the equality case of the Fenchel-type theorem.
\end{remark}
The paper is organized as follows. Section 2 reviews the necessary background on Sasakian sub-Riemannian 3-manifolds and sub-Riemannian 3-space forms, as well as Legendrian curves and geodesics therein. In Section 3, we establish two fundamental properties of the flow (\ref{1.2}). Section 4 derives the evolution equations for geometric quantities with respect to the flow (\ref{1.2}). In Section 5, we establish higher order curvature estimates for $\gamma(\cdot,t)$, and prove the long-time existence of the flow (\ref{1.2}). Section 6 analyzes the asymptotic behavior of the flow (\ref{1.2}) in Sasakian 3-manifolds and proves Theorems \ref{thm1.3} and \ref{thm1.4}. Finally, Section 7 investigates the length-preserving flow (\ref{1.3}) in $\mathbb{M}(0)$ and establishes the geometric identities in Propositions \ref{prop1.7} and \ref{prop1.8}.

%%%
%%% Section 2
%%%

\section{Preliminaries}

\noindent In this section we collect some facts about Sasakian sub-Riemannian 3-manifolds and the model spaces in sub-Riemannian 3-space forms. We also review some properties of Legendrian curves and geodesics in these settings. See \cite{B2002,DG2006} for more details. 

\subsection{ Sasakian sub-Riemannian 3-manifold}
\
\vglue-10pt
 \indent

Let $M$ be a contact 3-manifold with the contact 1-form $\theta$, and denote the contact distribution by $\mathcal{H}: =  \rm{Ker}(\theta) \subset TM$. The choice of a scalar product $g$ on $\mathcal{H}$ defines a sub-Riemannian structure on $M$. In this case we say that $(M,\theta, g)$ is a contact sub-Riemannian manifold. The 2-form $d \theta$ is given by
\begin{align*}
    d \theta(U, V)=U(\theta(V))-V(\theta(U))-\theta([U, V]),
\end{align*}
where $[U, V]$ is the Lie bracket of two vector fields on $M$. We always choose the orientation of $M$ induced by the nowhere vanishing 3-form $\theta \land  d\theta $.

The Reeb vector field associated to $\theta$ is the smooth vector field $T$ transversal to $\mathcal{H}$ determined by 
\begin{align*}
    \theta(T)=1, \qquad d \theta(T, U)=0.
\end{align*}
Let $J$ be the orientation-preserving $90^{\circ}$ rotation on the contact distribution $\mathcal{H}$. This defines a complex structure on $\mathcal{H}$ compatible with $\theta$, i.e., $J^2=-\mathrm{Id}$ and equality
\begin{align*}
    g(U, V)= \frac{1}{2}d \theta(U, J(V)), \qquad (U,V) \in  \mathcal{H}.
\end{align*}
Thus, $g(J(U), J(V))  = g(U, V)$, and $g$ is a Hermitian metric. In this case, $(M,\theta, -J)$ is a pseudohermitian strictly pseudo-convex CR-manifold as defined in \cite[p.~76]{B2002}. 

The canonical extension of $g$ is a Riemannian metric $g_{\theta}$ on $M$ by declaring $T$ to a unit vector orthogonal to $\mathcal{H}$. Moreover, the complex structure $J$ is extend to the whole tangent space $TM$ by $JT := 0$. Clearly, the Riemannian metric $g_{\theta}$ on $M$ takes the form
\begin{align*}
    g_{\theta}:=\frac{1}{2}d \theta(\cdot, J(\cdot))+\theta\otimes  \theta.
\end{align*}
For any vector fields $U$ and $V$ on $TM$, it is easy to check that
\begin{align}\label{2.1}
    &J^{2}(U) = -U + \theta(U)T,  \\
    &g_{\theta}(U,T)=\theta (U), \notag \\
    &g_{\theta}(J(U),J(V)) = g_{\theta}(U,V) - \theta(U)\theta(V), \notag \\
    &g_{\theta}(J(U),V) + g_{\theta}(U,J(V))=0 \notag.
\end{align}
In \cite{T1989}, Tanno introduced the Tanaka connection on Sasakian manifolds and proved the following lemma.
\begin{lemma}\cite{B2002,T1989}\label{lem2.1}
    The Tanaka connection $\nabla$ on a contact metric manifold $(M,\theta,-J,g_{\theta})$ is the unique linear connection such that
    \begin{align*}
        &\text{(i) } \nabla \theta = 0, \\
        &\text{(ii) } \nabla T = 0, \\
        &\text{(iii) } \nabla g_{\theta} = 0, \\
        &\text{(iv) } T_{\nabla}(U, V) = d\theta(U, V)\, T, \quad \text{for all } U, V \in \mathcal{H}, \\
        &\text{(v) } T_{\nabla}(T, JU) = -J\, T_{\nabla}(T, U), \quad \text{for all } U \in TM,
    \end{align*}
    where $T_{\nabla}$ denotes the torsion tensor of the Tanaka connection.
\end{lemma}
The tensor $\mathcal{Q}$ and the torsion operator $\sigma$ associated with $\nabla$ are defined respectively as
\begin{align*}
    &\mathcal{Q}(U, V) := (\nabla_V J)U, \quad \forall U, V \in TM, \\
    &\sigma (U) := T_{\nabla}(T, U). \quad \forall U \in TM.
\end{align*}
A Sasakian sub-Riemannian 3-manifold is a contact sub-Riemannian 3-manifold for which $\mathcal{Q} = 0$ and $\sigma = 0$. Moreover, combining the properties (iv) and (v) of Lemma \ref{lem2.1}, we have  
\begin{lemma}\cite{DG2006}\label{lem2.2}
    Let $T_{\nabla}$ denote the torsion tensor of Sasakian manifold $(M, \theta, g_{\theta}, -J)$ with respect to the Tanaka connection $\nabla$. Then,
    \begin{align}\label{2.2}
        T_{\nabla}(U, V)=-2 g_\theta(U, J V) T, \quad \text{for any } U,V \in TM.
    \end{align}
\end{lemma}

The Levi-Civita connection $\nabla^\theta$ of $\left(M, g_\theta\right)$ satisfies the equality $\nabla^\theta_{U} T = J(U)$, for any vector field $U$. Furthermore, by \cite[p.~86,91]{B2002} and the second equality in (\ref{2.1}), we get
\begin{align*}
    \nabla^\theta_{U} (J(V)) = J(\nabla^\theta_{U} V) + g_{\theta}(V,T)U -  g_{\theta}(U,V)T.
\end{align*}
This implies that the CR-structure induced by $J$ is integrable. From the Koszul's formulas in \cite[p.~38]{DG2006}, we get the relationship between $\nabla$ and  $\nabla^\theta$, 
\begin{align}\label{2.3}
    \nabla^\theta=\nabla-\frac{1}{2} d \theta \otimes T+\theta \odot J
\end{align}
where $\theta \odot J(U, V)=\frac{1}{2}(\theta(U) J V+\theta(Y) J V)$, for any $U, V \in T M$.

The curvature tensor of a Sasakian sub-Riemannian 3-manifold  $\left(M, \theta , g_\theta \right)$ associated with the Tanaka connection $\nabla$ is defined by 
\begin{align*}
    R(U,V)Z = \nabla_{V}\nabla_{U} Z - \nabla_{U}\nabla_{V}Z +\nabla_{[U,V]}Z
\end{align*} 
for all $U,V,Z \in TM $. The Webster scalar curvature $W$ is the most important invariant of Sasakian manifolds. Up to a constant, $W$ coincides with the sectional curvature of the contact plane with respect to $\nabla$. Consider an orthonormal frame $\left\{e_1, e_2=J e_1, T\right\} $ with respect to $g_\theta$, the Webster scalar curvature $W$ is defined by
\begin{align}\label{2.4}
    2W:= g_{\theta}\left( R(e_{1},e_{2})e_{2}, e_{1}\right).
\end{align}

Given a curve $\gamma : I \subset \mathbb{R} \to M$ and let $s$ be the arc-length parameter of $\gamma$. Then,  $\dot{\gamma}=\frac{d \gamma}{d s}$ is a unit tangent vector field along $\gamma$ with respect to $g_\theta$. The curvature $k$ and the contact normality $\tau$ of $\gamma$ are 
\begin{align}\label{2.5}
    k =g_\theta\left(\nabla_{\dot{\gamma}} \dot{\gamma}, J \dot{\gamma}\right), \qquad \tau =g_\theta(\dot{\gamma}, T).
\end{align}
If $\tau=0$ along $\gamma$, then $\gamma$ is a Legendrian curve, i.e., $\gamma$ is tangent to $\mathcal{H}$ everywhere. Let $e_1:=\dot{\gamma}$ and $e_2 =J e_1=J\dot{\gamma}$. Then, $\left\{\dot{\gamma}, J\dot{\gamma}, T\right\}$ forms an orientable orthonormal basis of $M$ along $\gamma$. The curvature of $\gamma$ can be expressed by
\begin{align}\label{2.6}
    \nabla_{\dot{\gamma}} \dot{\gamma}=k J \dot{\gamma}.
\end{align}
\begin{lemma}\label{lem2.3}
    Let $M$ be a Sasakian sub-Riemannian 3-manifold and $\gamma: I \to  M$ be a Legendrian curve parameterized by arc-length. Then, $\gamma$ is geodesic if and only if there is $\lambda \in \mathbb{R}$ such that $\gamma$ satisfying the equation
    \begin{align}\label{2.7}
        \nabla_{\dot{\gamma}}\dot{\gamma} + 2\lambda J\dot{\gamma} = 0.
    \end{align}
\end{lemma}
It is observed that all geodesics in the Sasakian sub-Riemannian 3-manifold have constant curvature. A variation of $\gamma$ is a $C^2$ map $F: I \times I^{\prime} \rightarrow M$, where $I^{\prime}$ is an open interval containing the origin, such that $F(s,0) = \gamma(s)$.  It was proved in \cite{C2012} that the vector field associated to the variation of geodesic satisfies a certain differential equation.

\begin{lemma}\label{lem2.4}
    Let $\gamma_{\varepsilon}(s) := F(s,\varepsilon)$ and $V_{\varepsilon}(s)$ be the vector field along $\gamma_{\varepsilon}$ given by $\partial F/ \partial \varepsilon(\varepsilon,s)$. Suppose $\gamma_{\varepsilon}(s)$ is the geodesic of curvature $k$ with the initial conditions $\gamma_{0}(s) = \gamma$ and $V_{0}(s) = \dot{\gamma_{0}}(s) $. Then $V_{\varepsilon}(s)$ satisfies
    \begin{align*}
        \frac{\partial^{3}}{\partial s^{3}} g_{\theta}(V_{\varepsilon}, T)+(k^{2}+ 2W)\frac{\partial}{ \partial s}g_{\theta}(V_{\varepsilon}, T)=0
    \end{align*}
    where $W$ is the Webster scalar curvature of $M$.
\end{lemma}

Moreover, either Lemmas 1.5 and 2.1 of \cite{SP2009} or Lemma 3.4 of \cite{C2012} provides an explicit computation of $g_{\theta}(V_{\varepsilon}, T)$ in sub-Riemannian 3-space forms.

\begin{lemma}\label{lem2.5} 
    Consider a geodesic $\gamma(s)$ of curvature $k$ in a Sasakian 3-manifold with constant Webster scalar curvature. Denote $V$ be the vector field corresponding to a variation of $\gamma$, and let $\alpha := k^2 + 2W$. Then, we have\\
    (i) For $\alpha < 0$,
    \begin{align*}
        g_{\theta}(V, T)(s) = \frac{1}{\alpha} (1-\cosh (\sqrt{-\alpha s})).
    \end{align*}
    (ii) For $\alpha = 0$,
    \begin{align*}
        g_{\theta}(V, T)(s) =\frac{1}{2} s^2.
    \end{align*}
    (iii) For $\alpha > 0$,
    \begin{align*}
        g_{\theta}(V, T)(s) = \frac{1}{\alpha} (1-\cos (\sqrt{\alpha s})).
    \end{align*}
\end{lemma}

\subsection{ Sub-Riemannian 3-space forms}
\
\vglue-10pt
\indent

The model spaces of sub-Riemannian 3-space forms $\mathbb{M}(K)$ are recalled briefly, with further details can be found in Section 7.4 of \cite{B2002}.

Let $\mathbb{N}(K)$ be the complete, simply connected, Riemannian surface of constant sectional curvature $K$. Explicitly,
\begin{align*}
    \mathbb{N}(1) = \mathbb{S}^2,\qquad \mathbb{N}(0) = \mathbb{R}^{2},\qquad \mathbb{N}(-1) = \{ (x,y) \in \mathbb{R}^{2} \mid x^{2} + y^{2} < 1\},
\end{align*}
where $\mathbb{S}^{2}$ is endowed with the standard sphere metric, $\mathbb{R}^{2}$ with the metric $d x^2+d y^2$, and $\mathbb{N}(-1)$ with the metric $\rho^2\left(d x^2+d y^2\right)$, where $\rho(x, y):=\left(1-\left(x^2+y^2\right)\right)^{-1}$.

The sub-Riemannian 3-space forms $\mathbb{M}(K)$ may be understood as sub-Riemannian counterparts to the Riemannian 2-space forms. Consequently, we define $\mathbb{M}(K)$ as follows. For $K =1$, let $\mathbb{M}(1) =  \mathbb{S}^{3} \subset \mathbb{R}^{4}$ equipped with the contact 1-form
\begin{align*}
    \theta = x_1 d y_1-y_1 d x_1+x_2 d y_2-y_2 d x_2, \qquad \{x_{1}, y_{1}, x_{2},y_{2}\} \in \mathbb{R}^{4}.
\end{align*}
For $K =-1,0$, let $\mathbb{M}(K) = \mathbb{N}(K) \times \mathbb{R}$ with the contact 1-form defined by
\begin{equation*}
    \theta = 
    \begin{cases}
        \rho(x d y-y d x)+d z, \qquad K =-1,  \\
        x d y-y d x +d z, \qquad K =0,
    \end{cases}
\end{equation*}
where $(x, y, z)$ are Euclidean coordinates in $\mathbb{R}^3$.

Let $\{ X_{1},X_{2},T\}$ be a basis on $\mathbb{M}(0)$, where $T$ is the Reeb vector field. We introduce in $\mathbb{M}(K)$ the sub-Riemannian metric $g$ on $\mathcal{H}$ for which $\{X_{1},X_{2} \}$ is a positively oriented orthonormal basis at each point. The associated complex structure $J$ satisfies $J(X_{1})=X_{2}$ and $J(X_{2})=-X_{1}$.  Similarly, $g$ and $J$ can be extend on $\mathbb{M}(K)$ and the whole tangent space $T\mathbb{M}(K)$, respectively.

Notice $\mathbb{M}(0)$ is the first Heisenberg group and the left-invariant vector field $\{ X_{1},X_{2},T\}$ is given by 
\begin{align*}
    X_{1}=\partial_x+y \partial_z, \quad X_{2}=\partial_y-x \partial_z, \quad T=\partial_z.
\end{align*}
By computing the Lie bracket of $X_{1}$ and $X_{2}$, we easily verify that
\begin{align*}
    -2T=[X_{1},X_{2}] = X_{1}X_{2}-X_{2}X_{1} .
\end{align*}
The contact distribution $\mathcal{H} = Span\{X_{1}, X_{2}\}$ and $T=\partial_z$ is the Reeb vector field. Moreover, we can express a curve $\gamma: S^{1} \to \mathbb{M}(0)$ in its coordinates functions, i.e.,
\begin{align*}
    \gamma=\gamma^1 \partial_x+\gamma^2 \partial_y+\gamma^3 \partial_z=\gamma^1 X_1+\gamma^2 X_2+\gamma^3 T ,
\end{align*}
and 
\begin{align*}
    \gamma^{\prime}(u)=\gamma_u^1 \partial_x+\gamma_u^2 \partial_y+\gamma_u^3 \partial_z =\gamma_u^1 X_1+\gamma_u^2 X_2+\left(\gamma_u^3+\gamma^1 \gamma_u^2-\gamma^2 \gamma_u^1\right) T.
\end{align*}
If $\gamma$ is a Legendrian curve in $\mathbb{M}(0)$, then $\gamma^{\prime}(u) \in \mathcal{H}$, which equivalent to 
\begin{align*}
    \gamma_u^3+\gamma^1 \gamma_u^2-\gamma^2 \gamma_u^1=0.
\end{align*}
The following results concerning geodesics in $\mathbb{M}(K)$ hold for $K =1,0$.  
\begin{lemma}[\cite{CDPT2007,AC2008,C2012}]\label{lem2.6}
    Let $\gamma$ be a complete geodesic of curvature $k$ in $\mathbb{M}(K)$. Then, $\gamma$ is either injective curves defined on $\mathbb{R}$, or closed circles with the same length. In particular,\\
    (i) In $\mathbb{M}(1)$, depending on an arithmetic condition involving $k$, $\gamma$ is either a closed circle of fixed length or it parameterizes a dense subset of a Clifford torus.\\
    (ii) In $\mathbb{M}(0)$, $\gamma$ is a horizontal straight line for $k=0$ and  a Euclidean helix with vertical axis for $k \neq 0$.
\end{lemma}

\section{Geometric Properties of the flow (\ref{1.2})}

\noindent In this section, we show that the flow (\ref{1.2}) increases the length of the evolving curve,  which reflects the expanding nature of this flow. Furthermore, we verify that the flow (\ref{1.2}) preserves the Legendrian condition.

\begin{proposition}\label{prop3.1}
    Along the flow (\ref{1.2}), the length of evolving curve $\gamma(\cdot, t)$ increases exponentially.
\end{proposition}
\noindent {\bf Proof}: Let $u$ be a parameter of $S^{1}$, which is independent of $t$. We define the arc-length parameter $s$ by  $s(u,t) = |\gamma^{\prime}(u,t)|$, and 
\begin{align}\label{3.1}
    \frac{\partial }{\partial s} = \frac{1}{|\gamma^{\prime}(u)|}\frac{\partial }{\partial u} , \qquad ds = |\gamma^{\prime}(u)| du.
\end{align}
By the definition of length, we have 
\begin{align*}
    \frac{d}{dt} \mathcal{L}_{t} =  \frac{d}{dt} \int_{\gamma_{t}} ds = \frac{d}{dt} \int_{S^{1}} |\gamma^{\prime}(u,t)\ du =  \int_{S^{1}} \frac{d}{dt}|\gamma^{\prime}(u,t)| du .
\end{align*}
In the following, we prove the evolution equation of $|\gamma^{\prime}(u,t)|$ along the flow (\ref{1.2}). Using (\ref{2.2}), we obtain
\begin{align*}
    \frac{\partial}{\partial t}|\gamma^{\prime}(u,t)|^2=&2 g_\theta\left(\nabla_t \gamma^{\prime}(u), \gamma^{\prime}(u)\right)=2 g_\theta\left(\nabla_u\left(\frac{\partial \gamma}{\partial t}\right)+T_{\nabla}\left(\frac{\partial \gamma}{\partial t}, \frac{\partial \gamma}{\partial u}\right), \gamma^{\prime}(u)\right)\\
    =& 2|\gamma^{\prime}(u)|^2g_\theta\left(\nabla_s\left(-k^{-1}J \dot{\gamma} - 2\int_{0}^{s} k^{-1}ds T 
    \right), \dot{\gamma}\right)\\
    &-4 |\gamma^{\prime}(u)|^2g_\theta\left(\frac{\partial \gamma }{\partial t} , J (\dot{\gamma})\right) g_\theta(T, \dot{\gamma}) \\
    =& 2|\gamma^{\prime}(u)|^2g_\theta\left(-\frac{\partial k^{-1}}{\partial s} J \dot{\gamma} - k^{-1} \nabla_{s} J \dot{\gamma} , \dot{\gamma}\right) \\
    =& 2|\gamma^{\prime}(u)|^2g_\theta\left(-\frac{\partial k^{-1}}{\partial s} J \dot{\gamma} - k^{-1} J(\nabla_{\dot{\gamma}} \dot{\gamma}), \dot{\gamma}\right)
\end{align*}
where the last equality relies on the fact that $\mathcal{Q} \equiv 0$ on Sasakian manifold $M$, i.e., $\nabla J = 0$. Furthermore, by the fourth equality in (\ref{2.1}) and (\ref{2.5}), we have
\begin{align*}
    g_{\theta}(J \dot{\gamma}, \dot{\gamma}) =  0, \qquad g_{\theta}(J(\nabla_{\dot{\gamma}} \dot{\gamma}), \dot{\gamma}) = - g_{\theta}(\nabla_{\dot{\gamma}} \dot{\gamma}, J \dot{\gamma})  = -k.
\end{align*}
Thus
\begin{align*}
    \frac{\partial}{\partial t}|\gamma^{\prime}(u,t)|^2= - 2|\gamma^{\prime}(u)|^2 k^{-1}g_{\theta}(J(\nabla_{\dot{\gamma}} \dot{\gamma}), \dot{\gamma}) = 2|\gamma^{\prime}(u,t)|^2,
\end{align*}
i.e.,
\begin{align}\label{3.2}
    \frac{\partial}{\partial t}|\gamma^{\prime}(u,t)| =  |\gamma^{\prime}(u,t)|.
\end{align}
Hence
\begin{align*}
    \frac{d}{dt} \mathcal{L}_{t} =  \int_{\gamma} |\gamma^{\prime}(u,t)| du =  \mathcal{L}_{t}.
\end{align*}
This completes the proof.
\hfill${\square}$

\begin{proposition}\label{prop3.2}
    If the initial curve $\gamma_{0}$ is Legendrian, then along the flow (\ref{1.2}), the evolving curve $\gamma_{t}$ is Legendrian for all $t>0$.
\end{proposition}
\noindent {\bf Proof}:  From the definition of Legendrian curves, we have $\tau(\cdot, 0) =g_\theta(\dot{\gamma}(\cdot,0), T) \equiv 0$. Along the flow (\ref{1.2}), we obtain
\begin{align*}
    \frac{\partial}{\partial t}\tau = g_{\theta}(\nabla_{t} \dot{\gamma}, T) 
\end{align*}
and 
\begin{align}\label{3.3}
    \nabla_{t} \dot{\gamma} & =\nabla_t\left(\frac{1}{\left|\gamma^{\prime}(u)\right|} \gamma^{\prime}(u)\right) \notag \\
    & =\frac{1}{\left|\gamma^{\prime}(u)\right|} \nabla_t \gamma^{\prime}(u) -\left(\frac{1}{\left|\gamma^{\prime}(u)\right|^2} \frac{\partial}{\partial t}\left|\gamma^{\prime}(u)\right|\right) \gamma^{\prime}(u) \notag \\
    & =\frac{1}{\left|\gamma^{\prime}(u)\right|} \nabla_u\left(\frac{\partial \gamma}{\partial t}\right)+\frac{1}{\left|\gamma^{\prime}(u)\right|} T_{\nabla}\left(\frac{\partial \gamma}{\partial t}, \frac{\partial \gamma}{\partial u}\right)  -\dot{\gamma} \notag \\
    &= \nabla_{s}\left(-k^{-1}J \dot{\gamma} 
    - 2\int_{0}^{s} k^{-1}ds T \right) + T_{\nabla}\left(\frac{\partial \gamma}{\partial t}, \dot{\gamma}\right) -\dot{\gamma} \notag \\
    & = - \frac{\partial k^{-1}}{\partial s}J \dot{\gamma} - k^{-1} J(\nabla_{\dot{\gamma}}\dot{\gamma}) -\dot{\gamma}.
\end{align}
Thus
\begin{align}\label{3.4}
    \frac{\partial}{\partial t}\tau = - \frac{\partial k^{-1}}{\partial s}g_{\theta}( J \dot{\gamma}, T) - k^{-1}   g_{\theta}(J(\nabla_{\dot{\gamma}}\dot{\gamma}), T) - g_{\theta} (\dot{\gamma}, T).
\end{align}
By the fourth equality in (\ref{2.1}), we get
\begin{align*}
    g_{\theta}( J\dot{\gamma}, T) =  - g_{\theta}( \dot{\gamma}, J(T)) = 0,\qquad g_{\theta}( J(\nabla_{\dot{\gamma}}\dot{\gamma}), T) =  - g_{\theta}( \nabla_{\dot{\gamma}}\dot{\gamma}, J(T)) = 0.
\end{align*}
Substituting the above equality into (\ref{3.4}) yields
\begin{align*}
    \frac{\partial}{\partial t}\tau = -  \tau.
\end{align*}
Then, by the initial condition $\tau (\cdot, 0)\equiv 0$, we get 
\begin{align*}
    \tau(\cdot, t)  = e^{-t}  \tau(\cdot, 0) = 0
\end{align*}
which implies that $\gamma_{t}$ is Legendrian for all $t>0$.
\hfill${\square}$

\begin{remark}\label{rem3.3}
    Since the solution  $\gamma(u,t)$ of the initial value problem (\ref{1.2}) is a Legendrian curve, we have $\nabla_{\dot{\gamma}} \dot{\gamma} = k J \dot{\gamma} $ holds for all $t>0$. It is not difficult to observe that (\ref{1.2}) is equivalent to 
\begin{equation*}
    \begin{cases}
        \frac{\partial \gamma}{\partial t} = -\frac{\nabla_{\dot{\gamma}}\dot{\gamma}}{|\nabla_{\dot{\gamma}}\dot{\gamma}|^{2}}
      - 2\int_{0}^{s} k^{-1}ds T,   \\
    \gamma(\cdot ,0) = \gamma_{0}.
    \end{cases}
\end{equation*}
In addition, from (\ref{2.3}), we have
\begin{align*}
    \nabla_{\dot{\gamma}} \dot{\gamma}=\nabla_{\dot{\gamma}}^\theta \dot{\gamma}-\theta(\dot{\gamma}) J (\dot{\gamma}) =  \nabla_{\dot{\gamma}}^\theta \dot{\gamma}.
\end{align*}
Hence, the flow (\ref{1.2}) differs from the one-dimensional case of the IMCF in $(M,g_{\theta})$ only by some lower-order terms. By the standard PDE theory, it follows that the flow (\ref{1.2}) admits a smooth solution for a short time.
\end{remark}

\section{Evolution equation }

\begin{lemma}\label{lem4.1}
    Along the flow (\ref{1.2}), we have the following evolution equations:
   \begin{align}
    &\nabla_{t} \dot{\gamma} = k^{-2} \frac{\partial k}{\partial s} J\dot{\gamma}. \label{4.1}\\
    &\nabla_{t}\nabla_{s} = \nabla_{s}\nabla_{t} - \nabla_{s} -2k^{-1}WJ . \label{4.2}\\
    &\frac{\partial }{\partial t}k^{2} = k^{-2} \frac{\partial^{2} k^{2}}{\partial s^{2}} - \frac{3}{2}k^{-4}\left(\frac{\partial k^{2}}{\partial s}\right)^{2} - 2k^{2}-4W. \label{4.3}\\
    &\frac{\partial }{\partial t}k^{-2} = k^{-2} \partial_{s}^{2}k^{-2} - \frac{1}{2}(\partial_{s}k^{-2})^{2}+ 2k^{-2}+4k^{-4}W.\label{4.4} 
   \end{align}
\end{lemma}
\noindent {\bf Proof}: From Proposition \ref{prop3.2}, $\nabla_{\dot{\gamma}}\dot{\gamma} = k J\dot{\gamma}$ and $\dot{\gamma} \in \mathcal{H}$. Since $ J^{2}\Big|_{\mathcal{H}} = -Id$, we get 
\begin{align}\label{4.5}
    J(\nabla_{\dot{\gamma}}\dot{\gamma}) = -k \dot{\gamma}.
\end{align}
Combining (\ref{3.3}) and (\ref{4.5}),  we obtain 
\begin{align*}
    \nabla_{t} \dot{\gamma} =- \frac{\partial k^{-1}}{\partial s} J \dot{\gamma} = k^{-2} \frac{\partial k}{\partial s} J \dot{\gamma}. 
\end{align*}

Now, we prove the relationship of covariant differentiations with respect to the arc-length $s$ and time $t$. A direct computation yields
\begin{align*}
    \nabla_{t}\nabla_{s} &= \nabla_{t}\left(\frac{1}{\left|\gamma^{\prime}(u)\right|}\nabla_{u}\right) =  \frac{1}{\left|\gamma^{\prime}(u)\right|}\nabla_{t}\nabla_{u} -  \nabla_{s} \\
    & = \nabla_{s}\nabla_{t} + R\left(\frac{\partial \gamma}{\partial t}, \dot{\gamma}\right) -  \nabla_{s}\\
\end{align*}
and 
\begin{align*}
    R\left(\frac{\partial \gamma}{\partial t}, \dot{\gamma}\right) = -k^{-1}R\left(J\dot{\gamma}, \dot{\gamma}\right) -2\int_{0}^{s} k^{-1} ds R\left(T, \dot{\gamma}\right) .
\end{align*}
Using the definition of $R$ and (\ref{2.4}), we get 
\begin{align*}
    &R\left(J\dot{\gamma}, \dot{\gamma}\right)e_{1} = g_{\theta}\left( R(e_{1},e_{2})e_{2}, e_{1}\right) e_{2} = 2W e_{2},\\
    &R\left(J\dot{\gamma}, \dot{\gamma}\right)e_{2} = -g_{\theta}\left( R(e_{1},e_{2})e_{2}, e_{1}\right) e_{1} = -2We_{1},\\
    &R\left(J\dot{\gamma}, \dot{\gamma}\right)T = 0.
\end{align*}
Thus
\begin{align*}
    R\left(J\dot{\gamma}, \dot{\gamma}\right) = 2W J.
\end{align*}
Similarly, we have $ R\left(T, \dot{\gamma}\right) = 0$. Using the above equalities, we obtain
\begin{align*}
    \nabla_{t}\nabla_{s} = \nabla_{s}\nabla_{t} - \nabla_{s} -2k^{-1}WJ.
\end{align*}

For notational convenience, we set $\partial_{s} : = \frac{\partial }{\partial s}$ and $\partial_{t} : = \frac{\partial }{\partial t}$. By (\ref{4.1}), (\ref{4.2}) and (\ref{4.5}), we obtain
\begin{align*}
    \frac{\partial }{\partial t}k^{2} &= 2g_{\theta}\left(\nabla_{t}\nabla_{s} \dot{\gamma}, \nabla_{\dot{\gamma}}\dot{\gamma}\right) \\
    &=  2g_{\theta}\left(\nabla_{s}\nabla_{t} \dot{\gamma}, \nabla_{\dot{\gamma}}\dot{\gamma} \right) - 2g_{\theta}(\nabla_{\dot{\gamma}}\dot{\gamma},\nabla_{\dot{\gamma}}\dot{\gamma}) - 4k^{-1}W g_{\theta}(J\dot{\gamma}, \nabla_{\dot{\gamma}}\dot{\gamma})\\
    & = 2g_{\theta}\left(\nabla_{s}\left( k^{-2}\partial_{s} kJ\dot{\gamma} \right), \nabla_{\dot{\gamma}}\dot{\gamma}\right) - 2k^{2} -4W \\
    & = 2k^{-1}\partial_{s}^{2}k + 2k \partial_{s}k^{-2} \partial_{s}k - 2k^{2}-4W\\
    & = 2k^{-1}\partial_{s}^{2}k  - 4k^{-2}(\partial_{s}k)^{2}- 2k^{2}-4W.
\end{align*}
Since $\partial_{s}k^{2} = 2k\partial_{s} k$ and $\partial_{s}^{2}k^{2} = 2k \partial_{s}^{2}k + 2(\partial_{s}k)^{2}$, equation (\ref{4.3}) follows immediately. Then, applying (\ref{4.3}), a direct computation yields (\ref{4.4}).
\hfill${\square}$

\begin{lemma}\label{lem4.2}
    The evolution of the higher derivative of $k^{-2}$ is given by 
    \begin{align}\label{4.6}
        \frac{\partial }{\partial t}\left| \partial_{s}^{m} k^{-2}\right|^2 =& k^{-2}\frac{\partial^{2}}{\partial s^{2}}\left| \partial_{s}^{m} k^{-2}\right|^2- 2k^{-2}\left| \partial_{s}^{m+1} k^{-2}\right|^{2} \\
        &+ 2(2-m)\left| \partial_{s}^{m} k^{-2}\right|^2 + 2\partial_{s}^{m} k^{-2} \sum^{m}_{l=1}\binom{m}{l} \partial_{s}^{l} k^{-2} \partial_{s}^{m-l+2} k^{-2}\notag \\
        &- \partial_{s}^{m} k^{-2}\sum^{m}_{l=0}\binom{m}{l}\partial_{s}^{l+1} k^{-2} \partial_{s}^{m-l+1} k^{-2} \notag \\
        &+8\partial_{s}^{m} k^{-2}\sum^{m}_{i=0}\binom{m}{i} \partial_{s}^{i} W\sum^{m-i}_{j=0}\binom{m-i}{j}\partial_{s}^{j} k^{-2} \partial_{s}^{m-i-j} k^{-2}, \notag
    \end{align} 
    where $\partial_{s}^{i} W = \frac{\nabla^{i}}{\partial s^{i}} W(\gamma(s))$ and $\partial_{s} W = g_{\theta}\left(\nabla W, \frac{\partial \gamma}{\partial s}\right)$ .      
\end{lemma}
\noindent {\bf Proof}: From (\ref{3.1}) and (\ref{3.2}), we get
\begin{align*}
    \frac{\partial}{\partial t}\frac{\partial}{\partial s} = \frac{\partial}{\partial t} \left(\frac{1}{\left|\gamma^{\prime}(u,t)\right|} \frac{\partial}{\partial u}\right) =  \frac{\partial}{\partial s}\frac{\partial}{\partial t} - \frac{\partial}{\partial s}.
\end{align*}
Hence
\begin{align*}
    \partial_{t}\partial_{s}^{m}k^{-2} = \partial_{s}\partial_{t}\partial_{s}^{m-1}k^{-2} - \partial_{s}\partial_{s}^{m-1}k^{-2} = \cdots = \partial_{s}^{m}\partial_{t}k^{-2} - m \partial_{s}^{m}k^{-2}
\end{align*}
and 
\begin{align*}
    \frac{\partial }{\partial t} \left| \partial_{s}^{m} k^{-2}\right|^2 = 2\partial_{s}^{m} k^{-2} \partial_{s}^{m} \left(\partial_{t} k^{-2}\right) - 2m \left| \partial_{s}^{m} k^{-2}\right|^{2}.
\end{align*}
Using (\ref{4.4}), we have 
\begin{align*}
    \partial_{s}^{m} \left(\partial_{t}k^{-2}\right)  =& \sum^{m}_{l=0}\binom{m}{l} \partial_{s}^{l} k^{-2} \partial_{s}^{m-l+2} k^{-2} - \frac{1}{2}\sum^{m}_{l=0}\binom{m}{l}\partial_{s}^{l+1} k^{-2} \partial_{s}^{m-l+1} k^{-2}  \\
    &+ 2 \partial_{s}^{m} k^{-2}+ 4\sum^{m}_{i=0}\binom{m}{i} \partial_{s}^{i} W\sum^{m-i}_{j=0}\binom{m-i}{j}\partial_{s}^{j} k^{-2} \partial_{s}^{m-i-j} k^{-2}
\end{align*}
and 
\begin{align}\label{4.7}
    \frac{\partial }{\partial t}\left| \partial_{s}^{m} k^{-2}\right|^2 =& 2\partial_{s}^{m} k^{-2} \sum^{m}_{l=0}\binom{m}{l} \partial_{s}^{l} k^{-2} \partial_{s}^{m-l+2} k^{-2} + 2(2-m)\left| \partial_{s}^{m} k^{-2}\right|^2\\
    &- \partial_{s}^{m} k^{-2} \sum^{m}_{l=0}\binom{m}{l}\partial_{s}^{l+1} k^{-2} \partial_{s}^{m-l+1} k^{-2}   \notag\\
    & + 8\partial_{s}^{m} k^{-2}\sum^{m}_{i=0}\binom{m}{i} \partial_{s}^{i} W\sum^{m-i}_{j=0}\binom{m-i}{j}\partial_{s}^{j} k^{-2} \partial_{s}^{m-i-j} k^{-2}.  \notag
\end{align}
Furthermore, 
\begin{align*}
    \frac{\partial^{2}}{\partial s^{2}}\left| \partial_{s}^{m} k^{-2}\right|^{2} = 2 \partial_{s}^{m+2} k^{-2}\partial_{s}^{m} k^{-2} + 2\left| \partial_{s}^{m+1} k^{-2}\right|^{2},
\end{align*}
which is implicitly contained within the $l=0$ case of the first term on the right-hand side of (\ref{4.7}). Thus, (\ref{4.6}) is obtained.
\hfill${\square}$

\section{Long-time existence }

\noindent  In this section, we establish the higher order curvature estimates for $\gamma_{t}$ by constructing suitable auxiliary functions for the cases $W>0$, $W<0$, and $W=0$. These estimates allow us to prove the long-time existence of the flow (\ref{1.2}) in Sasakian 3-manifolds. We first show that $k^{2}$ is uniformly bounded.

\begin{lemma}\label{lem5.1}
    Let $\gamma_{t}$ be a solution to the flow (\ref{1.2}) in Sasakian 3-manifold $M$. Then, there exists a constant $c_{0} > 0$ depending on $\gamma_{0}$ and $w_{0} : =\max_{M}|W|$ such that 
    \begin{align}\label{5.1}
        \max_{\gamma_{t}}k^{2} \leq c_{0}, \qquad t\in [0,T_{0}).
    \end{align}
\end{lemma}
\noindent {\bf Proof}: Let $F = e^{2t}\left(k^{2} - 2w_{0}\right)$, by (\ref{4.3}), we have 
\begin{align*}
    \frac{\partial }{\partial t} F  &= k^{-2}\partial_{s}^{2} F- \frac{3}{2}e^{-2t}k^{-4}\left(\partial_{s} F\right)^{2} - 4e^{2t}(W+w_{0})\\
    &\leq k^{-2}\partial_{s}^{2} F- \frac{3}{2}e^{-2t}k^{-4}\left(\partial_{s} F\right)^{2}.
\end{align*}
Denote $F_{\max}(t): = \max_{\gamma_{t}} F(\cdot, t)$, we have
\begin{align*}
    \frac{\partial }{\partial t} F_{\max} \leq 0.
\end{align*}
Thus
\begin{align*}
    e^{2t} \left(\max_{\gamma_{t}}k^{2}  - 2w_{0}\right) = F_{\max}(t) \leq F_{\max}(0) =  \left(\max_{\gamma_{0}}k^{2}  - 2w_{0}\right),
\end{align*}
i.e.,
\begin{align*}
    \max_{\gamma_{t}}k^{2} \leq e^{-2t}\left(\max_{\gamma_{0}}k^{2}- 2w_{0}\right)+ 2w_{0}.
\end{align*}
We choose  $c_{0}:=2w_{0}$ if $\max_{\gamma_{0}}k^{2}- 2w_{0} \leq 0$. Otherwise, we set $c_{0} :=\max_{\gamma_{0}}k^{2}$.
\hfill${\square}$

Next, we establish the higher order estimates of $k$ in Sasakian 3-manifolds with positive Webster scalar curvature. For notational convenience, we denote
\begin{align*}
   w_{l} = \max_{M} \left|\nabla^{i} W\right|^{2}.
\end{align*}
\begin{theorem}\label{thm5.2}
    Suppose that the Webster scalar curvature $W$ of $M$ is bounded below by a positive constant $W_{0}$, i.e., $W \geq W_{0} > 0$. Then, for each $m \geq 1$, there exists a constant $c_{m}$ depending on $c_{0}$, $w_{0}, \cdots, w_{m}$ and $W_{0}$, such that 
    \begin{align*}
        \max_{\gamma_{t}} \left| \partial_{s}^{m} k\right|^{2} \leq c_{m},\qquad t\in [0,T_{0}).
    \end{align*}
\end{theorem}
\noindent {\bf Proof}: 
(1) When $m =1$, we construct the function 
\begin{align*}
    \phi = \left|\partial_{s} k^{-2}\right|^{2} + A_{1}k^{-4},
\end{align*}
where $A_{1}$ is a constant to be determined later. Using (\ref{4.4}), we obtain
\begin{align}\label{5.2}
    \frac{\partial }{\partial t}k^{-4} = k^{-2}\partial_{s}^{2} k^{-4} - 3k^{-2}\left(\partial_{s} k^{-2}\right)^{2} + 4k^{-4} + 8Wk^{-6}.
\end{align}
From Lemma \ref{lem4.2},  we have 
\begin{align}
    \frac{\partial }{\partial t}\left|\partial_{s} k^{-2}\right|^{2} =& k^{-2}\partial_{s}^{2}\left|\partial_{s} k^{-2}\right|^{2} -2k^{-2}\left|\partial_{s}^{2} k^{-2}\right|^{2} + 2\left|\partial_{s} k^{-2}\right|^{2}  \label{5.3} \\
    &+ 8k^{-4} \partial_{s}W \partial_{s}k^{-2}+16Wk^{-2}\left|\partial_{s} k^{-2}\right|^{2}  \notag \\
    \leq &  k^{-2}\partial_{s}^{2}\left|\partial_{s} k^{-2}\right|^{2} -2k^{-2}\left|\partial_{s}^{2} k^{-2}\right|^{2}  \label{5.4} \\
    &+ \left(2k^{2} + 16W + 16(\partial_{s}W)^{2}\right) k^{-2}\left|\partial_{s} k^{-2}\right|^{2} + k^{-6} \notag
\end{align}
Thus, by (\ref{5.2}) and (\ref{5.3}), we have 
\begin{align*}
    \frac{\partial }{\partial t}\phi = &k^{-2} \partial_{s}^{2}\phi - 2k^{-2}|\partial_{s}^{2}k^{-2}|^{2} + \left(2 + 16Wk^{-2} - 3A_{1}k^{-2}\right) \left|\partial_{s} k^{-2}\right|^{2} \\
    &+ 8k^{-4}\partial_{s}W \partial_{s}k^{-2} + 4A_{1}k^{-4} +8A_{1}Wk^{-6}.
\end{align*}
Denote $\phi_{\max} := \max_{\gamma_{t}} \phi(\cdot,t)$, then
\begin{align*}
    \partial_{s} \phi_{\max} = 2\partial_{s}k^{-2}\partial_{s}^{2}k^{-2} + 2A_{1}k^{-2}\partial_{s}k^{-2} =  0 ,
\end{align*}
which implies that either $\partial_{s}k^{-2} = 0$ or $\partial_{s}^{2}k^{-2} = - A_{1}k^{-2}$ holds.

If $\partial_{s}k^{-2} = 0$ at the spatial maximum point of $\phi$, we choose $A_{1} = 0$. Then, $\frac{\partial }{\partial t}\phi_{\max} \leq 0$, i.e., $\phi_{\max}(t) \leq \phi_{\max}(0)$. Hence
\begin{align*}
   \max_{\gamma_{t}} \left|\partial_{s} k^{-2}\right|^{2} \leq \max_{\gamma_{0}} \left|\partial_{s} k^{-2}\right|^{2}.
\end{align*}

If $\partial_{s}^{2}k^{-2} = - A_{1}k^{-2}$ at the spatial maximum point of $\phi$, using (\ref{5.2}) and (\ref{5.4}), we have 
\begin{align*}
    \frac{\partial }{\partial t}\phi \leq  &k^{-2} \partial_{s}^{2}\phi + \left(2k^{2} + 16W + 16(\partial_{s}W)^{2}- 3A_{1}\right) k^{-2}\left|\partial_{s} k^{-2}\right|^{2} \\
    &+ \left(1 + 4A_{1}k^{2} +  8A_{1}W - 2A_{1}^{2}\right)k^{-6} \\
    \leq & k^{-2} \partial_{s}^{2}\phi + \left(2c_{0} + 16w_{0} + 16w_{1}- 3A_{1}\right) k^{-2}\left|\partial_{s} k^{-2}\right|^{2} \\
    &+ \left(1 + 4A_{1}c_{0} +  8A_{1}w_{0} - 2A_{1}^{2}\right)k^{-6} \\
\end{align*}
Choosing
\begin{align*}
    A_{1} := \max \left\{ \frac{1}{3}(2c_{0} + 16w_{0} + 16w_{1}) , 2c_{0} + 4w_{0} + (8w_{0})^{-1}\right\},
\end{align*}
we have 
\begin{align*}
    2c_{0} + 16w_{0} + 16w_{1}- 3A_{1} \leq 0, 
\end{align*} 
and 
\begin{align*}
    1+4A_{1}c_{0} +  8A_{1}w_{0}-2A_{1}^{2} \leq -\frac{c_{0}}{2w_{0}} - \frac{1}{32w_{0}^{2}} <0.
\end{align*}
Thus, $\partial_{t} \phi_{\max} \leq 0$ and 
\begin{align*}
    \max_{\gamma_{t}} \left(\left|\partial_{s} k^{-2}\right|^{2} + A_{1}k^{-4}\right)\leq \max_{\gamma_{0} }\left(\left|\partial_{s} k^{-2}\right|^{2} + A_{1}k^{-4}\right) : = c_{1}.
\end{align*}
In summary, as $c_{1} > \max_{\gamma_{0} } \left|\partial_{s} k^{-2}\right|^{2}$, it follows that
\begin{align}\label{5.5}
    \max_{\gamma_{t}} \left|\partial_{s} k^{-2}\right|^{2} \leq c_{1}.
\end{align}

(2) When $m \geq 2$, we need to construct an auxiliary function to eliminate bad terms. Let 
\begin{align*}
    \bar{\phi} :=  \left|\partial_{s} k^{-2}\right|^{2} -W_{0}^{-1} k^{-4}  - \alpha_{1}k^{-2},
\end{align*}
where $\alpha_{1}$ is a positive constant to be determined later. Combining (\ref{4.4}), (\ref{5.2}) and (\ref{5.3}), we obtain
\begin{align}\label{5.6}
    \partial_{t} \bar{\phi} =& k^{-2}\partial_{s}^{2}\bar{\phi} - 2k^{-2}\left|\partial_{s}^{2} k^{-2}\right|^{2} + 2\left|\partial_{s} k^{-2}\right|^{2}  + 8k^{-4} \partial_{s}W \partial_{s}k^{-2}+16Wk^{-2}\left|\partial_{s} k^{-2}\right|^{2} \notag \\
    &+ 3W_{0}^{-1}k^{-2}\left|\partial_{s} k^{-2}\right|^{2} - 4W_{0}^{-1}k^{-4} - 8\frac{W}{W_{0}}k^{-6} \notag \\
    &+\frac{1}{2} \alpha_{1}\left|\partial_{s} k^{-2}\right|^{2} - 2\alpha_{1}k^{-2} - 4\alpha_{1}Wk^{-4} \notag \\
    \leq& k^{-2}\partial_{s}^{2}\bar{\phi} - 2k^{-2}\left|\partial_{s}^{2} k^{-2}\right|^{2} + \left( 2k^{2} + 16(\partial_{s}W)^{2} +16 W + 3W_{0}^{-1}\right) k^{-2}\left|\partial_{s} k^{-2}\right|^{2} \\
    &-7k^{-6} - 2\alpha_{1}k^{-2} + \frac{1}{2} \alpha_{1}c_{1} \notag .
\end{align}
Denote $\beta := 2c_{0} + 16w_{1} +16w_{0} + 3W_{0}^{-1}$, we have $\beta > 0$ and (\ref{5.6}) can be reduced to 
\begin{align*}
     \partial_{t} \bar{\phi} \leq& k^{-2}\partial_{s}^{2}\bar{\phi} - 2k^{-2}\left|\partial_{s}^{2} k^{-2}\right|^{2} +\beta k^{-2} |\partial_{s} k^{-2}|^{2} -7k^{-6} - 2\alpha_{1}k^{-2} + \frac{1}{2} \alpha_{1}c_{1} \\
     =&k^{-2}\partial_{s}^{2}\bar{\phi} - 2k^{-2}\left|\partial_{s}^{2} k^{-2}\right|^{2}  -7k^{-6} + \left(\beta |\partial_{s} k^{-2}|^{2}  - 2\alpha_{1}\right)k^{-2}+ \frac{1}{2} \alpha_{1}c_{1}.
\end{align*}
By choosing $\alpha_{1} = \beta c_{1}$ and using (\ref{5.5}), the above inequality is equivalent to 
\begin{align}\label{5.7}
    \partial_{t} \bar{\phi} 
    \leq k^{-2}\partial_{s}^{2}\bar{\phi} - 2k^{-2}\left|\partial_{s}^{2} k^{-2}\right|^{2}-7k^{-6}  - \beta k^{-2} \left|\partial_{s} k^{-2}\right|^{2}+ \frac{1}{2} \beta c_{1}^{2}.
\end{align}

In the following, we construct the function
\begin{align*}
    \phi_{2} = \left|\partial_{s}^{2} k^{-2}\right|^{2} +A_{2}\left|\partial_{s} k^{-2}\right|^{2} + \alpha_{2} \bar{\phi}
\end{align*}
where $A_{2}$ and $\alpha_{2}$ are positive constants to be determined later. From Lemma \ref{lem4.2}, we have 
\begin{align}\label{5.8}
    \frac{\partial }{\partial t}\left|\partial_{s}^{2} k^{-2}\right|^{2} =& k^{-2}\partial_{s}^{2}\left|\partial_{s}^{2} k^{-2}\right|^{2} -2k^{-2}\left|\partial_{s}^{3} k^{-2}\right|^{2} + 2(\partial_{s} k^{-2}) (\partial_{s}^{2} k^{-2}) (\partial_{s}^{3} k^{-2})  \notag \\
    &+ 16Wk^{-2}|\partial_{s}^{2} k^{-2}|^{2} + 16W|\partial_{s} k^{-2}|^{2}(\partial^{2}_{s} k^{-2}) \notag\\
    &+ 32(\partial_{s}W)k^{-2}(\partial_{s} k^{-2})(\partial^{2}_{s} k^{-2}) +8(\partial^{2}_{s}W)k^{-4}(\partial^{2}_{s} k^{-2}) \notag \\
    \leq & k^{-2}\partial_{s}^{2}\left|\partial_{s}^{2} k^{-2}\right|^{2} -k^{-2}\left|\partial_{s}^{3} k^{-2}\right|^{2}  + k^{2}|\partial_{s} k^{-2}|^{2} |\partial^{2}_{s} k^{-2}|^{2} + 16Wk^{-2}|\partial^{2}_{s} k^{-2}|^{2} \\
    &+ 16W^{2}|\partial_{s} k^{-2}|^{2}|\partial^{2}_{s} k^{-2}|^{2} +4|\partial_{s} k^{-2}|^{2} + 16(\partial_{s}W)^{2}k^{-2}|\partial^{2}_{s} k^{-2}|^{2} \notag \\
    &+ 16k^{-2}|\partial_{s} k^{-2}|^{2} + 16(\partial^{2}_{s} W)^{2} k^{-2}|\partial^{2}_{s} k^{-2}|^{2} + k^{-6}\notag.
\end{align}
Combining (\ref{5.4}), (\ref{5.7}) and the above inequality, we get 
\begin{align}\label{5.9}
    \partial_{t} \phi_{2} \leq& k^{-2}\partial_{s}^{2}\phi_{2} - k^{-2}\left|\partial_{s}^{3} k^{-2}\right|^{2} + \left[k^{4}|\partial_{s} k^{-2}|^{2} + 16W + 16W^{2}k^{2}|\partial_{s} k^{-2}|^{2} \right. \\
    &\left. + 16 (\partial_{s} W)^{2} + 16(\partial^{2}_{s} W)^{2} - 2A_{2} -2\alpha_{2}\right]k^{-2}|\partial^{2}_{s} k^{-2}|^{2} \notag \\
    &+ \left(4k^{2} +16 +2k^{2}A_{2} +  16A_{2}W + 16A_{2}(\partial_{s}W)^{2}  - \beta \alpha_{2}\right)k^{-2}|\partial_{s} k^{-2}|^{2} \notag\\
    &+(1+A_{2} -7\alpha_{2})k^{-6} + \frac{1}{2}\alpha_{2}\beta c_{1}^{2} \notag
\end{align}
Let $\beta_{2} := c_{0}^{2}c_{1} + 16w_{0} + 16w_{0}^{2}c_{0}c_{1} + 16w_{1} + 16w_{2}$. We choose $A_{2} = \beta_{2} $ and 
\begin{align*}
    \alpha_{2}:= \max \left\{  \beta^{-1}\left(4c_{0}+16+2c_{0}\beta_{2} +16w_{0}\beta_{2} + 16w_{1}\beta_{2}\right), \frac{1}{7}(\beta_{2}+1)\right\}.
\end{align*}
Hence, (\ref{5.9}) can be reduced to
\begin{align}\label{5.10}
    \partial_{t} \phi_{2} \leq& k^{-2}\partial_{s}^{2}\phi_{2} - k^{-2}\left|\partial_{s}^{3} k^{-2}\right|^{2}  - \beta_{2}k^{-2}|\partial^{2}_{s} k^{-2}|^{2} + \frac{1}{2}\alpha_{2}\beta c_{1}^{2} \notag \\
    \leq & k^{-2}\partial_{s}^{2}\phi_{2} - k^{-2}\left|\partial_{s}^{3} k^{-2}\right|^{2}  - \beta_{2}c_{0}^{-1}|\partial^{2}_{s} k^{-2}|^{2} +   \frac{1}{2}\alpha_{2}\beta c_{1}^{2}
\end{align}
where the last inequality is obtained by applying (\ref{5.1}). 

From $\phi_{2} = \left|\partial_{s}^{2} k^{-2}\right|^{2} +\beta_{2}\left|\partial_{s} k^{-2}\right|^{2} + \alpha_{2} \bar{\phi}$ and  $\bar{\phi} =  \left|\partial_{s} k^{-2}\right|^{2} -W_{0}^{-1} k^{-4}  - \beta c_{1}k^{-2}$, we have 
\begin{align}\label{5.11}
    \left|\partial_{s}^{2} k^{-2}\right|^{2} &= \phi_{2} - (\beta_{2}+\alpha_{2})|\partial_{s} k^{-2}|^{2} + \alpha_{2}W_{0}^{-1}k^{-4} + \alpha_{2}\beta c_{1}k^{-2}\\
    &\geq \phi_{2} - (\beta_{2}+\alpha_{2})|\partial_{s} k^{-2}|^{2}. \notag
\end{align}
Thus
\begin{align}\label{5.12}
    -\beta_{2}c_{0}^{-1}\left|\partial_{s}^{2} k^{-2}\right|^{2} \leq  -\beta_{2}c_{0}^{-1} \phi_{2} + (\beta_{2}+\alpha_{2})\beta_{2}c_{0}^{-1}c_{1}
\end{align}
Therefore, substituting (\ref{5.12}) into (\ref{5.10}) yields
\begin{align*}
    \partial_{t} \phi_{2} \leq& k^{-2}\partial_{s}^{2}\phi_{2} - k^{-2}\left|\partial_{s}^{3} k^{-2}\right|^{2}  - \beta_{2}c_{0}^{-1}\phi_{2} + \varepsilon_{1}
\end{align*}
where $\varepsilon_{1} = \frac{1}{2}\alpha_{2}\beta c_{1}^{2}  + (\beta_{2}+\alpha_{2})\beta_{2}c_{0}^{-1}c_{1}$. By the maximum principal, we get
\begin{align*}
    \max_{\gamma_{t}} \phi_{2}(s,t) \leq (\max_{\gamma_{0}} \phi_{2} - \varepsilon_{1} c_{0} \beta_{2}^{-1})e^{-\beta_{2}c_{0}^{-1}t} + \varepsilon_{1} c_{0} \beta_{2}^{-1} \leq  \max_{\gamma_{0}} \phi_{2} + \varepsilon_{1} c_{0} \beta_{2}^{-1}
\end{align*}
Denote $\Lambda_{2} : =  \max_{\gamma_{0}} \phi_{2} + \varepsilon_{1} c_{0} \beta_{2}^{-1}$, and from (\ref{5.11}), we have 
\begin{align}\label{5.13}
    k^{4}\left|\partial_{s}^{2} k^{-2}\right|^{2} +k^{4}(\beta_{2} +\alpha_{2})\left|\partial_{s} k^{-2}\right|^{2} &\leq \Lambda_{2} k^{4} + \alpha_{2}W_{0}^{-1} + \alpha_{2}\beta c_{1}k^{2}  \notag \\
    &\leq  \Lambda_{2} c_{0}^{2} + \alpha_{2}W_{0}^{-1} + \alpha_{2}\beta c_{1}c_{0} : =\Lambda_{2}^{*}.
\end{align}
Since $\partial_{s}^{2} k^{-2} = 6k^{-4} (\partial_{s}k)^{2} - 2k^{-3}\partial_{s}^{2}k$, we have 
\begin{align*}
    \left|\partial_{s}^{2} k^{-2}\right|^{2} &= 36k^{-8} (\partial_{s}k)^{4} + 4k^{-6}(\partial^{2}_{s}k)^{2} - 24k^{-7}(\partial_{s}k)^{2}\partial^{2}_{s} k \\
    &\geq  3k^{-6}(\partial^{2}_{s}k)^{2}  - 108k^{-8}(\partial_{s}k)^{4}.
\end{align*}
Thus, using (\ref{5.1}), (\ref{5.5}) and (\ref{5.13}), we get  
\begin{align*}
    k^{-2}|\partial^{2}_{s}k|^{2} \leq \frac{1}{3}k^{4}\left|\partial_{s}^{2} k^{-2}\right|^{2} + 36k^{-4}(\partial_{s}k)^{4} \leq \frac{1}{3}\Lambda_{2}^{*} + \frac{9}{4}c_{0}^{4}c_{1}^{2}:=c_{2}^{*}
\end{align*}
and 
\begin{align}\label{5.14}
    |\partial^{2}_{s}k|^{2} \leq c_{2}^{*} k^{2} \leq c_{2}^{*} c_{0}:= c_{2}.
\end{align}

(3) When $m\geq 3$, we construct the function
\begin{align*}
    \phi_{m} = |\partial^{m}_{s}k^{-2}|^{2} + \sum_{i=3}^{m}A_{i} |\partial^{i-1}_{s}k^{-2}|^{2}  + \alpha_{m} \bar{\phi}
\end{align*}
where $A_{i}$ and $\alpha_{m}$ are positive constants. From Lemma \ref{4.2}, we have
\begin{align*}
    \frac{\partial }{\partial t}\left| \partial_{s}^{m} k^{-2}\right|^{2}=&k^{-2}\frac{\partial^{2}}{\partial s^{2}}\left(\left| \partial_{s}^{m} k^{-2}\right|^{2}\right) - 2k^{-2}\left| \partial_{s}^{m+1} k^{-2}\right|^{2} + 2(2-m)\left| \partial_s^{m}k^{-2}\right|^2\\
    &+(2m-2)\partial_{s}k^{-2}\partial_{s}^{m}k^{-2} \partial_{s}^{m+1}k^{-2}+2\sum_{l = 2}^{m} \binom{m}{l} \partial_{s}^{l}k^{-2} \partial_{s}^{m-l+2}k^{-2} \partial_{s}^{m}k^{-2}\\
    & - \partial_{s}^{m}k^{-2} \sum_{l = 1}^{m-1} \binom{m}{l}\partial_{s}^{l+1}k^{-2}\partial_{s}^{m-l+1}k^{-2}  \\
    &+ 8 \partial_{s}^{m} k^{-2} \sum^{m}_{i=0}\binom{m}{i} \partial_{s}^{i} W\sum^{m-i}_{j=0}\binom{m-i}{j}\partial_{s}^{j} k^{-2} \partial_{s}^{m-i-j} k^{-2}.
\end{align*}
Since
\begin{align*}
    (2m-2)\partial_{s}k^{-2}\partial_{s}^{m}k^{-2} \partial_{s}^{m+1}k^{-2} \leq (m-1)^{2}k^{2}|\partial_{s}k^{-2}|^{2}\left| \partial_{s}^{m} k^{-2}\right|^{2} + k^{-2}\left| \partial_{s}^{m+1} k^{-2}\right|^{2},
\end{align*}
we get 
\begin{align}\label{5.15}
    \frac{\partial }{\partial t}\left| \partial_{s}^{m} k^{-2}\right|^{2}
    \leq &k^{-2}\frac{\partial^{2}}{\partial s^{2}}\left(\left| \partial_{s}^{m} k^{-2}\right|^{2}\right) - k^{-2}\left| \partial_{s}^{m+1} k^{-2}\right|^{2} \\
    &+\left[2(2-m) + (m-1)^{2}k^{2}|\partial_{s}k^{-2}|^{2}\right]\left| \partial_s^{m}k^{-2}\right|^2 \notag \\
    &+2 \partial_{s}^{m}k^{-2} \sum_{l = 2}^{m} \binom{m}{l} \partial_{s}^{l}k^{-2} \partial_{s}^{m-l+2}k^{-2} \notag \\
    &- \partial_{s}^{m}k^{-2} \sum_{l = 1}^{m-1} \binom{m}{l}\partial_{s}^{l+1}k^{-2}\partial_{s}^{m-l+1}k^{-2}  \notag \\
    &+ \partial_{s}^{m} k^{-2} 8\sum^{m}_{i=0}\binom{m}{i} \partial_{s}^{i} W\sum^{m-i}_{j=0}\binom{m-i}{j}\partial_{s}^{j} k^{-2} \partial_{s}^{m-i-j} k^{-2}.\notag
\end{align}
Suppose we have bounds 
\begin{align*}
     k^{4} |\partial_{s}^{i} k^{-2}|^{2} \leq \Lambda^{*}_{i}, \qquad t \in [0,T_{0})
\end{align*}
for $2\leq i \leq m-1$. Then, this assumption and (\ref{5.5}) allow us to find positive constants $A_{i}$ and $\alpha_{m}$ to enable
\begin{align*} 
    \frac{\partial }{\partial t}\phi_{m} \leq k^{-2}\partial_{s}^{2} \phi_{m} -  \varepsilon_{m}\phi_{m} + \varepsilon_{m+1},
\end{align*}
where $\varepsilon_{m} $ and $\varepsilon_{m+1}$ are positive constants which depending on $\{w_{i}\}_{i=0}^{m}$, $\{ \Lambda_{i}^{*}\}_{i=2}^{m-1}$, $c_{0}$, $c_{1}$ and $W_{0}$. By the maximum principle, we obtain
\begin{align*}
    \max_{\gamma_{t}}\phi_{m} \leq \max_{\gamma_{0}} \phi_{m} + \frac{\varepsilon_{m+1}}{\varepsilon_{m}}:=\Lambda_{m}. 
\end{align*}
Consequently,
\begin{align*}
    k^{4}|\partial_{s}^{m} k^{-2}|^{2} \leq k^{4} \phi_{m} +\alpha_{m}(W_{0}^{-1} + \beta c_{1}k^{2})  \leq c_{0}^{2}\Lambda_{m} + \alpha_{m}(W_{0}^{-1} + \beta c_{1}c_{0}) : = \Lambda_{m}^{*}.
\end{align*}

Moreover, from $\partial_{s}^{2} k^{-2} = 6k^{-4} (\partial_{s}k)^{2} - 2k^{-3}\partial_{s}^{2}k$, we deduce that $\partial_{s}^{m} k^{-2}$ can be written as $\partial_{s}^{m} k^{-2} =  f(k^{-1},\partial_{s}^{i} k) - 2k^{-3}\partial_{s}^{m}k$, where $f(k^{-1}, \partial_{s}^{i} k)$ is a polynomial in $k^{-1}$ and derivatives $\partial_{s}^{i} k$ of order less than $m$. Hence
\begin{align*}
    |\partial_{s}^{m} k^{-2}|^{2} \geq 3k^{-6}|\partial_{s}^{m} k|^{2} -3f^{2}
\end{align*}
and
\begin{align*}
    k^{-2}|\partial_{s}^{m} k|^{2} \leq \frac{1}{3}k^{4}|\partial_{s}^{m} k^{-2}|^{2} + k^{4}f^{2} \leq \frac{1}{3}\Lambda_{m}^{*}+c_{0}^{2}f^{2}
\end{align*} 
Similarly, from (\ref{5.5}) and the above assumption, it follows that $f^{2}$ is uniformly bounded. Therefore, there exists a constant $c^{*}_{m}$ depending on $\{\Lambda_{i}^{*}\}_{i=2}^{m}$, $c_{0}$ and $c_{1}$ such that 
\begin{align*}
    |\partial^{m}_{s}k|^{2}  \leq k^{2}c^{*}_{m} \leq c_{0}c^{*}_{m} :=c_{m}.
\end{align*}
This completes the proof of Theorem \ref{thm5.2}.
\hfill${\square}$

In the following, we focus on establishing higher order curvature estimates in Sasakian 3-manifolds with negative Webster scalar curvature. The first step is to demonstrate the uniform boundedness of $k^{-2}$.
\begin{lemma}\label{lem5.3}
    Suppose that the Webster scalar curvature $W$ of $M$ is bounded above by a negative constant $\bar{W}_{0}$, i.e., $W \leq \bar{W}_{0} < 0$. Then, there exists a constant $C_{0} >0$ depending on $\bar{W}_{0}$ and $\gamma_{0}$ such that
    \begin{align}\label{5.16}
        \max_{\gamma_{t}}k^{-2} \leq C_{0}, \qquad t \in [0, T_{0}).
    \end{align}
\end{lemma} 
\noindent {\bf Proof}: From (\ref{4.4}), we have 
\begin{align*}
    \partial_{t}k^{-2} =& k^{-2} \partial_{s}^{2}k^{-2} - \frac{1}{2}(\partial_{s}k^{-2})^{2}+ 2k^{-2}+4Wk^{-4} \\
    \leq & k^{-2} \partial_{s}^{2}k^{-2}  - \left(2|\bar{W}_{0}|^{\frac{1}{2}}k^{-2} - |\bar{W}_{0}|^{-\frac{1}{2}}\right)^{2} -  2k^{-2}  + |\bar{W}_{0}|^{-1} \\
    \leq &  k^{-2} \partial_{s}^{2}k^{-2} -  2k^{-2}  + |\bar{W}_{0}|^{-1} 
\end{align*}
Thus, by the maximum principal, we get 
\begin{align}\label{5.17}
    \max_{\gamma_{t} }k^{-2} \leq \left(\max_{\gamma_{0}}k^{-2} - \frac{1}{2|\bar{W}_{0}|}\right)e^{-2t} + \frac{1}{2|\bar{W}_{0}|}.
\end{align}
Denote $C_{1} :=\max_{\gamma_{0}}k^{-2} +( 2|\bar{W}_{0}|)^{-1}$, then (\ref{5.16}) follows from (\ref{5.17}). 
\hfill${\square}$

\begin{theorem}\label{thm5.4}
    Under the assumption made in Lemma \ref{lem5.3}, 
    there exist a constant $C_{m}$ depending on $w_{0}, \cdots, w_{m}$, $C_{0}$ and $\bar{W}_{0}$, such that 
    \begin{align*}
        \max_{\gamma_{t}} \left| \partial_{s}^{m} k^{-2}\right|^{2} \leq C_{m},\qquad t\in [0,T_{0}).
    \end{align*}
    holds for each $ m \geq 1$.
\end{theorem}
\noindent {\bf Proof}: (1) When $m =1$, we construct the function
\begin{align*}
    \psi = \left|\partial_{s} k^{-2}\right|^{2} + B_{1}k^{-2},
\end{align*}
where $B_{1}$ is a constant to be determined later. Using (\ref{5.4}) and (\ref{5.16}), we find that there exists a positive constant $\delta_{1}$ such that 
\begin{align}\label{5.18}
    \frac{\partial }{\partial t}\left|\partial_{s} k^{-2}\right|^{2} \leq  k^{-2}\partial_{s}^{2}\left|\partial_{s} k^{-2}\right|^{2} -2k^{-2}\left|\partial_{s}^{2} k^{-2}\right|^{2} + \delta_{1}\left(\left|\partial_{s} k^{-2}\right|^{2} +1\right),
\end{align}
where $\delta_{1}$ depending on $w_{1}$, $C_{0}$. Combining (\ref{4.4}) and (\ref{5.18}), we have
\begin{align*}
    \frac{\partial }{\partial t} \psi &\leq k^{-2}\partial_{s}^{2}\psi +  \delta_{1} \left(\left|\partial_{s} k^{-2}\right|^{2} +1\right) - \frac{1}{2}B_{1} \left|\partial_{s} k^{-2}\right|^{2} + 2B_{1}k^{-2}+4B_{1}Wk^{-4} \\
    & = k^{-2}\partial_{s}^{2}\psi+\left( \delta_{1} - \frac{1}{2}B_{1}\right)\psi + \left(\frac{1}{2}B_{1}^{2}+2B_{1}-B_{1}\delta_{1}\right)k^{-2} + \delta_{1} + 4B_{1}Wk^{-4}
\end{align*}
Choosing $ B_{1} = 4\delta_{1} $, we obtain 
\begin{align*}
    \frac{\partial }{\partial t} \psi &\leq k^{-2}\partial_{s}^{2}\psi - \delta_{1} \psi + 4\delta_{1}(\delta_{1} +2)k^{-2} + \delta_{1}  \\
    & \leq k^{-2}\partial_{s}^{2}\psi  - \delta_{1} \psi +4\delta_{1}(\delta_{1} +2)C_{0} + \delta_{1}
\end{align*}

Let $\psi_{\max}: = \max_{\gamma_{t}} \psi$. By the maximum principle, we get 
\begin{align*}
    \psi_{\max}(t) &\leq e^{-\delta_{1} t}\left[\max_{\gamma_{0}}\psi - 4(\delta_{1} + 2)C_{0} -1\right]  +4(\delta_{1} + 2)C_{0} +1 \\
    &\leq \max_{\gamma_{0}}\psi +  4(\delta_{1} + 2)C_{0} +1 
\end{align*}
Thus, by Lemma \ref{lem5.3}, there exists a positive constant $C_{1}$ such that 
\begin{align}\label{5.19}
    \max_{\gamma_{t}} \left|\partial_{s} k^{-2}\right|^{2} \leq C_{1}, \qquad t\in [0,T_{0}),
\end{align}
where $C_{1}$ depending on $\gamma_{0}$, $w_{1}$ and $C_{0}$. 

(2) When $m \geq 2$, we construct the function
\begin{align*}
    \psi_{m} = \left| \partial_{s}^{m}k^{-2}\right|^2 + B_{m} \left| \partial_{s}^{m-1}k^{-2}\right|^2,
\end{align*} 
where $B_{m}$ is a positive constants to be determined later. Assume that there exist constants $C_{i}$ such that 
\begin{align*}
    \left| \partial_{s}^{i} k^{-2}\right|^{2} \leq C_{i},\qquad   t \in [0,T_{0})
\end{align*}
for all $i=1,\cdots,m-1$. Then, combining this assumption and (\ref{5.15}), we can find a constant $\delta_{m} >0$ such that  
\begin{align*}
    \frac{\partial }{\partial t}\left| \partial_{s}^{m} k^{-2}\right|^{2}\leq & k^{-2}\frac{\partial^{2}}{\partial s^{2}}\left| \partial_{s}^{m} k^{-2}\right|^{2} - k^{-2}\left| \partial_{s}^{m+1} k^{-2}\right|^{2} +\delta_{m}\left(\left| \partial_{s}^{m} k^{-2}\right|^{2} +1\right).
\end{align*}
Hence,
\begin{align*} 
    \frac{\partial }{\partial t}\psi_{m} \leq k^{-2}\partial_{s}^{2} \psi_{m} - k^{-2}\left| \partial_{s}^{m+1} k^{-2}\right|^{2} + (\delta_{m} - B_{m})\left| \partial_{s}^{m} k^{-2}\right|^{2} + \delta_{m+1} 
\end{align*}
where $\delta_{m+1}$ is a constant depending on $w_{0},\cdots,w_{m}$, $C_{0}$ and $\bar{W}_{0}$. Choosing $B_{m} = 2\delta_{m}$, we have 
\begin{align*}
    \frac{\partial }{\partial t}\psi_{m}  \leq k^{-2}\partial_{s}^{2} \psi_{m} - k^{-2}\left| \partial_{s}^{m+1} k^{-2}\right|^{2} -  \delta_{m} \psi_{m}+ 2\delta_{m}^{2}C_{m-1} +  \delta_{m+1}.
\end{align*}
Thus, by the maximum principal, we obtain that $\psi_{m}$ is uniformly bounded and 
\begin{align*}
    \max_{\gamma_{t}}\left| \partial_{s}^{m} k^{-2}\right|^{2} \leq C_{m}, \qquad t \in [0, T_{0})
\end{align*}
where constant $C_{m}$ depending on $w_{0},\cdots, w_{m}$ and $C_{0}$. 
\hfill${\square}$

Finally, we establish the higher order curvature estimates on Sasakian 3-manifolds with vanishing Webster scalar curvature.
\begin{theorem}\label{thm5.5}
    Let $M$ be a Sasakian 3-manifold with vanishing Webster scalar curvature. Then, for each $m \geq 1$, there exist a constant $\tilde{c}_{m}$ depending on $w_{0}, \cdots, w_{m}$ and $c_{0}$, such that 
    \begin{align*}
        \max_{\gamma_{t}} \left| \partial_{s}^{m} k^{-2}\right|^{2} \leq \tilde{c}_{m},\qquad t\in [0,T_{0}).
    \end{align*}
\end{theorem}
\noindent {\bf Proof}: (1) When $m =1$, we construct the function 
\begin{align*}
    \tilde{\phi} = \left|\partial_{s} k^{-2}\right|^{2} + D_{1} k^{-4},
\end{align*}
where $D_{1}$ is a constant to be determined later. Using (\ref{5.2}) and (\ref{5.3}), we have 
\begin{align*}
    \frac{\partial }{\partial t} \tilde{\phi}  = k^{-2} \partial_{s}^{2} \tilde{\phi} - 2k^{-2}|\partial_{s}^{2}k^{-2}|^{2} + \left(2k^{2}- 3D_{1}\right) k^{-2}\left|\partial_{s} k^{-2}\right|^{2}  + 4D_{1}k^{-4}.
\end{align*}
Denote $\tilde{\phi}_{\max} := \max_{\gamma_{t}} \tilde{\phi}(\cdot,t)$, as in the proof of Theorem \ref{thm5.2} for the case $m=1$, one can find a nonnegative constant $D_{1}$ such that $\partial_{t} \tilde{\phi}_{\max} \leq 0$, i.e., $\max_{\gamma_{t}} \tilde{\phi} \leq \max_{\gamma_{0}} \tilde{\phi} $. Hence, there exists a positive constant $\tilde{c}_{1}$ depending on $\gamma_{0}$ such that 
\begin{align*}
    \max_{\gamma_{t}} \left|\partial_{s} k^{-2}\right|^{2} \leq \tilde{c}_{1}, \qquad t \in [0,T_{0}).
\end{align*}

(2) When $m \geq 2$, we construct the function 
\begin{align*}
    \tilde{\phi}_{m} = \left|\partial^{m}_{s} k^{-2}\right|^{2} + D_{m}|\partial_{s}^{m-1} k^{-2}|^{2},
\end{align*}
where $D_{m}$ is a positive constant. Suppose we have bounds
\begin{align*}
    |\partial_{s}^{i} k^{-2}|^{2} \leq \tilde{c}_{i},\qquad  1\leq i \leq m-1.
\end{align*}
Then, combining this assumption and (\ref{5.1}), there exist suitable positive constants $\eta_{1}$ and $\eta_{2}$ such that 
\begin{align*}
    \frac{\partial }{\partial t}\tilde{\phi}_{m} \leq k^{-2} \partial_{s}^{2} \tilde{\phi}_{m} - (\eta _{1} -D_{m})k^{-2}|\partial_{s}^{m} k^{-2}|^{2} + \eta _{2}
\end{align*}
Choose $D_{m} = 2\eta_{1}$, we have 
\begin{align*}
    \frac{\partial }{\partial t}\tilde{\phi}_{m} \leq k^{-2} \partial_{s}^{2} \tilde{\phi}_{m} - \eta _{1}c^{-1}_{0} \tilde{\phi} + 2\eta _{1}^{2}c_{0}^{-1}\tilde{c}_{m-1}+ \eta _{2}.
\end{align*}
The maximum principal implies that $\tilde{\phi}_{m}$ and then also $\left|\partial^{m}_{s} k^{-2}\right|^{2}$ are uniformly bounded.
\hfill${\square}$

\noindent {\bf Proof of Theorem \ref{thm1.1}}: Using $\nabla_{s} \dot{\gamma} = k J\dot{\gamma}$ and $\nabla_{s}J\dot{\gamma} = - k \dot{\gamma}$, we have 
\begin{align}\label{5.20}
    \nabla^{2}_{s} \dot{\gamma} = \partial_{s}kJ\dot{\gamma} - k^{2} \dot{\gamma}.
\end{align}
Let $\nabla^{m}_{s} \dot{\gamma} :=  f(s,t)J\dot{\gamma} +g(s,t)\dot{\gamma}$. By repeatedly differentiating (\ref{5.20}) with respect to $s$, one obtains that $f(s,t)$ and $g(s,t)$ are polynomials in $k, \partial_{s}k,\cdots, \partial_{s}^{m-1}k$. This implies that, if $k^{2}$ and $|\partial_{s}^{m} k|^{2}$ are bounded, then $\left|\nabla_{s}^{m} \dot{\gamma}\right|^{2}$ must also be bounded.

For $W>0$, combining Lemma \ref{lem5.1} with Theorem \ref{thm5.2}, we can obtain that $\left|\nabla_{s}^{m} \dot{\gamma}\right|^{2}$ is uniformly bounded on $[0,T_{0})$ for every $m \geq 1$. 

For $W \leq 0$, Theorems \ref{thm5.4} and \ref{thm5.5} establish that $|\partial_{s}^{m} k^{-2}|^{2}$ is bounded. Since $\partial_{s}^{m} k$ can be expressed in terms of $k$ and $\partial_{s}^{i} k^{-2}$ for $i=1,\dots,m$, it follows from Lemma \ref{lem5.1} that $|\partial_{s}^{m} k|^{2}$ is uniformly bounded. This implies that $ \left|\nabla_{s}^{m} \dot{\gamma}\right|^{2}$ is uniformly bounded on $[0,T_{0})$.

Since $ \left|\nabla_{s}^{m} \dot{\gamma}\right|^{2}$ is uniformly bounded on $[0,T_{0})$, we can deduce that $\left| \frac{\partial }{\partial t}\nabla^{m}_{s} \dot{\gamma}\right|^{2}$ are also bounded. Therefore, at time $T_{0}$, the tangent vector $\dot{\gamma}$ has a well-defined limit, which by integration yields a smooth curve $\gamma$. Moreover, by the short-time existence theorem (see Remark \ref{3.3}), $\gamma$ can be extended to a smooth solution on an interval strictly larger than $[0, T_{0}]$. This completes the proof of Theorem \ref{thm1.1}.
\hfill${\square}$

\section{ Asymptotic behavior}

\indent  In this section, we analyze the asymptotic behavior of the flow (\ref{1.2}) on Sasakian sub-Riemannian 3-manifolds.

\noindent {\bf Proo of Theorem \ref{thm1.3}}: Since the Webster scalar curvature $W \geq 0$, from (\ref{4.3}), we have
\begin{align*}
    \partial_{t}k^{2} \leq k^{-2}\partial_{s}^{2} k^{2} - \frac{3}{2}k^{-4} |\partial_{s} k^{2}|^{2} - 2k^{2}.
\end{align*} 
Let $k^{2}_{max} := \max_{\gamma_{t}} k^{2}(s,t)$. Applying the maximum principle, we obtain
\begin{align*}
    k^{2}_{\max} \leq e^{-2t} \max_{\gamma_{0}} k^{2}.
\end{align*}
Thus
\begin{align}\label{6.1}
    \lim_{t \to \infty} k^{2} = 0.
\end{align}
In the case $W>0$, as shown in the proof of Theorem \ref{thm5.2}, there exists a positive constant $c^{*}_{l}$ such that 
\begin{align}\label{6.2}
    |\partial_{s}^{l} k|^{2} \leq c^{*}_{l} k^{2}, \qquad   t \in [0,\infty),
\end{align}
holds for $l \geq 1$. Hence, combining (\ref{6.1}) and (\ref{6.2}), we obtain
\begin{align*}
    \lim_{t \to \infty} |\partial_{s}^{l} k|^{2}  = 0,  \qquad \text { for every } l \geq 1.
\end{align*}

In the case $W = 0$, Theorem \ref{thm5.5} shows that $|\partial_{s}^{m} k^{-2}|^{2}$ is uniformly bounded. Since $\partial_{s}^{m} k$ can be written in terms of $k$ and $\partial_{s}^{i} k^{-2}$ for $i=1,\dots,m$, it follows from (\ref{6.1}) that $\lim_{t \to \infty} |\partial_{s}^{l} k|^{2} = 0$ holds for any $l \geq 1$. 

Therefore, for $W \geq 0$, $\lim_{t \to \infty} |\partial_{s}^{l} k|^{2} = 0$, which implies that $\lim_{t \to \infty} |\nabla_{s}^{l} \dot{\gamma}| = 0$. As $M$ is complete, there exist compact sets $\{ \mathcal{K}_{i}\}_{i=1}^{\infty}$ such that $M = \cup_{i=1}^{\infty}  \mathcal{K}_{i}$. Let $\gamma_{i}(\cdot,t) := \gamma (\cdot,t) \cap \mathcal{K}_{i}$. By the Arzel$\grave{\rm a}$-Ascoli Theorem, there exists a subsequence curves $\gamma_{i}(\cdot,t_{j}) \longrightarrow  \gamma_{\mathcal{K}}(\cdot, \infty)$. From Lemma \ref{lem2.3}, $\gamma_{\mathcal{K}}(\cdot, \infty)$ is a geodesic of vanishing curvature. 
\hfill${\square}$

\noindent {\bf Proof of Theorem \ref{thm1.4}}: Let $f(s,t) := k^{2}+2W$. Using the fact that $W$ is a negative constant and (\ref{4.3}), we have 
\begin{align*}
    \partial_{t} f  = k^{-2}\partial_{s}^{2}f - \frac{3}{2}k^{-4}|\partial_{s} f|^{2} - 2f.
\end{align*}
Thus, by maximum principal, we have 
$\max_{\gamma_{t}}f  \leq e^{-2t} \max_{\gamma_{0}} f$, i.e.,
\begin{align}\label{6.3}
    k^{2} \leq e^{-2t} \max_{\gamma_{0}} f -2W.
\end{align}

Let $g := k^{-2} - (2|W|)^{-1}$. From (\ref{4.4}), we have 
\begin{align*}
    \partial_{t} g= k^{-2}\partial_{s}^{2}g - \frac{1}{2}|\partial_{s} g|^{2} -(2|W|^{\frac{1}{2}}k^{-2} - |W|^{-\frac{1}{2}})^{2}- 2g \leq k^{-2}\partial_{s}^{2}g -2g. 
\end{align*}
Thus, the maximum principal implies that $\max_{\gamma_{t}}g \leq e^{-2t} \max_{\gamma_{0}} g$, i.e., 
\begin{align}\label{6.4}
    k^{-2} \leq e^{-2t} \max_{\gamma_{0}} g -(2W)^{-1}
\end{align}
Combining (\ref{6.3}) and (\ref{6.4}), we obtain
\begin{align*}
    (e^{-2t} \max_{\gamma_{0}} g -(2W)^{-1} )^{-1} \leq k^{2} \leq e^{-2t} \max_{\gamma_{0}} f -2W.
\end{align*}
Therefore,
\begin{align*}
    \lim_{t \to \infty} k^{2} = -2W
\end{align*}
Similarly, there exists a subsequence curves $\gamma_{i}(\cdot,t_{j}) \longrightarrow  \gamma_{\mathcal{K}}(\cdot, \infty)$ and $\gamma_{\mathcal{K}}(\cdot, \infty)$ is a geodesic of curvature $|k|=\sqrt{-2W}$. This completes the proof of Theorem \ref{thm1.4}.
\hfill${\square}$

\section{The flow (\ref{1.2}) in the first Heisenberg group and its applications}

\indent  In this section, we first introduce the relationship between the Legendrian curve expanding flow (\ref{1.2}) in $\mathbb{M}(0)$ and the IMCF of planar curves. Next, we obtain the length-preserving flow (\ref{1.3}) through a Heisenberg dilation of the flow (\ref{1.2}). Finally, using the geometric properties of the length-preserving flow (\ref{1.3}), we establish two geometric identities in $\mathbb{M}(0)$.

\subsection{ Relationship between the flow (\ref{1.2}) and the IMCF of planar curves }
\
\vglue-10pt
 \indent

\begin{lemma}\cite{PS2021}\label{lem7.1}
    Let $\gamma$ be a Legendrian curve in $\mathbb{M}(0)$ and the projection of $\gamma$ onto $x y$-plane denote by $\gamma^{*}$. Then 
    \begin{align*}
        &\left|\gamma^{\prime}\right|^2=\left|\left(\gamma^*\right)^{\prime}\right|_{\mathbb{R}^{2}}^2=\left(\gamma_u^1\right)^2+\left(\gamma_u^2\right)^2 ,\\
        &\dot{\gamma}=\gamma_s^1 X_{1}+\gamma_s^2 X_{2} ,\qquad
        \dot{\gamma}^*=\gamma_s^1 \partial_x+\gamma_s^2 \partial_y,\\
        &J \dot{\gamma}=-\gamma_s^2 X_{1}+\gamma_s^1 X_{2} ,\qquad \mathbf{n}^*=-\gamma_s^2 \partial_x+\gamma_s^1 \partial_y, \\
        &k=k^*=\gamma_s^1 \gamma_{s s}^2-\gamma_s^2 \gamma_{s s}^1.
    \end{align*}
\end{lemma}
Obviously, the curvature of $\gamma \subset \mathbb{M}(0)$ is the signed curvature of the plane curve $\gamma^{*} \subset \mathbb{R}^{2}$. Moreover, we have the following results:
\begin{lemma}\label{lem7.2}
    Let $ \gamma^{*} : S^{1} \times [0,\infty) \to \mathbb{R}^{2}$ be the projection of $\gamma(\cdot,t)$ onto $x y$-plane. If $\gamma(\cdot,t)$ is a solution to the flow (\ref{1.2}), then $\gamma^{*}$ evolves by the following flow
    \begin{align*}
        \frac{\partial \gamma^{*}}{\partial t} =  -(k^*)^{-1}\mathbf{n}^{*}.
    \end{align*}
\end{lemma}
\noindent {\bf Proof}: By direct computation, we have 
\begin{align*}
    \frac{\partial \gamma}{\partial t} & =\gamma_t^1 \partial_x+\gamma_t^2 \partial_y+\gamma_t^3 \partial_z \\
    & =\gamma_t^1 X_{1}+\gamma_t^2 X_{2}+\left(\gamma_t^3+\gamma^1 \gamma_t^2-\gamma^2 \gamma_t^1\right) T.
\end{align*} 
Since $\gamma(\cdot,t)$ is a solution to the flow (\ref{1.2}), it follows that 
\begin{align*}
    &g_{\theta}\left(\frac{\partial \gamma}{\partial t}, J\dot{\gamma}\right) = -\gamma_t^1\gamma_s^2+ \gamma_t^2\gamma_s^1 = - k^{-1},\qquad g_{\theta}\left(\frac{\partial \gamma}{\partial t}, \dot{\gamma}\right) = \gamma_t^1\gamma_s^1+ \gamma_t^2\gamma_s^2=0.   
\end{align*}
Thus
\begin{align}
    &\left\langle \frac{\partial \gamma^{*}}{\partial t}, \mathbf{n}^{*} \right\rangle _{\mathbb{R}^{2}}  = -\gamma_t^1\gamma_s^2+ \gamma_t^2\gamma_s^1 =g_{\theta}\left(\frac{\partial \gamma}{\partial t}, J\dot{\gamma}\right)= - k^{-1} = -(k^*)^{-1}, \label{7.1}\\
    &\left\langle \frac{\partial \gamma^{*}}{\partial t}, \dot{\gamma^{*}} \right\rangle _{\mathbb{R}^{2}} = \gamma_t^1\gamma_s^1+ \gamma_t^2\gamma_s^2 =g_{\theta}\left(\frac{\partial \gamma}{\partial t}, \dot{\gamma}\right) = 0 .\label{7.2}
\end{align}
Since $\{\dot{\gamma^{*}}, \mathbf{n}^{*}\}$ forms an orthonormal basis of $\mathbb{R}^{2}$ along $\gamma^{*}$, we obtain
\begin{align*}
    \frac{\partial \gamma^{*}}{\partial t} =  -(k^*)^{-1}\mathbf{n}^{*}.
\end{align*}
This completes the proof.
\hfill${\square}$

\begin{lemma}\label{lem7.3}
    If the initial curve $ \gamma_{0} : =\gamma (u,0) = (\gamma_{0}^{1} (u), \gamma_{0}^{2} (u), \gamma_{0}^{3} (u))$ and $ \gamma^{*}(u,t)$ is a solution of the following initial problem
    \begin{equation}\label{7.3}
        \begin{cases}
            \frac{\partial \gamma^{*}}{\partial t} =  -(k^*)^{-1}\mathbf{n}^{*}, \\
            \gamma^{*}(u,0) =  \gamma^{*}_{0} : = (\gamma_{0}^{1} (u), \gamma_{0}^{2} (u)).
        \end{cases}
    \end{equation}
    Then, there exists a function $f:S^{1}\times [0,\infty) \to \mathbb{R}$ with $f(u,0) = \gamma_{0}^{3} (u)$ such that $\gamma(u,t) = (\gamma^{*}, f(u,t))$ is a solution of the flow (\ref{1.2}) in $\mathbb{M}(0)$ with the initial value $\gamma_{0} $.
\end{lemma}
\noindent {\bf Proof}: Let 
\begin{align*}
    &f(u,t) =  f(0,t) + \int_{0}^{u} \left( \gamma_{w}^{1}\gamma^{2} - \gamma_{w}^{2}\gamma^{1}\right)(w,t) dw,\\
    &f(0,t) = \gamma_{0}^{3} (u) + \int_{0}^{t} \left(\gamma_{v}^{1}\gamma^{2} - \gamma_{v}^{2}\gamma^{1}\right)(0,v) d v.
\end{align*}
From Lemma \ref{lem7.1}, we have
\begin{align*}
    -\gamma_t^1\gamma_u^2+ \gamma_t^2\gamma_u^1 = -(k^*)^{-1} \left|\left(\gamma^*\right)^{\prime}(u,t)\right|_{\mathbb{R}^{2}}= - k^{-1} |\gamma^{\prime}(u,t)| .
\end{align*}
In addition,
\begin{align*}
    \partial_{t}\left( \gamma_{u}^{1}\gamma^{2} - \gamma_{u}^{2}\gamma^{1}\right) &=  \partial_{t}(\gamma_{u}^{1}) \gamma^{2} - \partial_{t}(\gamma_{u}^{2})\gamma^{1} + \left( \gamma_{u}^{1}\gamma_{t}^{2} - \gamma_{u}^{2}\gamma_{t}^{1}\right) \\
    &=  \partial_{u}(\gamma_{t}^{1} \gamma^{2} - \gamma_{t}^{2}\gamma^{1} )+ 2\left( \gamma_{u}^{1}\gamma_{t}^{2} - \gamma_{u}^{2}\gamma_{t}^{1}\right) .
\end{align*}
Thus,
\begin{align*}
    \frac{\partial }{\partial t}f(u,t) &= \frac{\partial }{\partial t}f(0,t) + \int_{0}^{u} \partial_{t}\left( \gamma_{w}^{1}\gamma^{2} - \gamma_{w}^{2}\gamma^{1}\right)(w,t) dw\\
    &=\left(\gamma_{t}^{1}\gamma^{2} - \gamma_{t}^{2}\gamma^{1}\right)(0,t) +\left(\gamma_{t}^{1}\gamma^{2} - \gamma_{t}^{2}\gamma^{1}\right)\bigg|_{(0,t)}^{(u,t)}  + 2\int_{0}^{u} \left( \gamma_{w}^{1}\gamma_{t}^{2} - \gamma_{w}^{2}\gamma_{t}^{1}\right) dw \\
    &=\left(\gamma_{t}^{1}\gamma^{2} - \gamma_{t}^{2}\gamma^{1}\right)(u,t) - 2\int_{0}^{u} k^{-1} |\gamma^{\prime}(w,t)| dw \\
    & = \left(\gamma_{t}^{1}\gamma^{2} - \gamma_{t}^{2}\gamma^{1}\right)(u,t) - 2\int_{0}^{s} k^{-1} ds ,
\end{align*}
and 
\begin{align*}
    \frac{\partial \gamma}{\partial t} &=\gamma_t^1 \partial_x+\gamma_t^2 \partial_y+\partial_{t}f \partial_z =\gamma_t^1 X_{1}+\gamma_t^2 X_{2}+\left(\partial_{t}f+\gamma^1 \gamma_t^2-\gamma^2 \gamma_t^1\right) T \\
    &=\gamma_t^1 X_{1}+\gamma_t^2 X_{2} - 2\left(\int_{0}^{s} k^{-1} ds\right) T.
\end{align*}
Since $\left\{\dot{\gamma}, J\dot{\gamma}, T\right\}$ forms an orthonormal basis of $\mathbb{M}(0)$ along $\gamma$, and combining (\ref{7.1}) and (\ref{7.2}), we have
\begin{align*}
    \frac{\partial \gamma}{\partial t} =  - k^{-1}J\dot{\gamma} - 2\left(\int_{0}^{s} k^{-1} ds\right) T.
\end{align*}
This completes the proof.
\hfill${\square}$

\begin{remark}\label{rem7.4}
    To specify a unique flow, we impose the condition $f(t,0) = f(0,0) =\gamma_{0}^{3} (0)$. Since the projection $\gamma^{*}(u,t)$ is the unique solution to the initial problem (\ref{7.3}), the components $\gamma^{1}(u,t)$ and $\gamma^{2}(u,t)$ are uniquely determined. The component $\gamma^{3}(u,t)$ is then uniquely determined by $\gamma^{3}_{u} = \gamma^{2}\gamma^{1}_{u} -\gamma^{1}\gamma^{2}_{u}$ together with $f(t,0) = \gamma_{0}^{3} (0)$. 
\end{remark}

\subsection{ Length-preserving flow (\ref{1.3}) and its applications}
\
\vglue-10pt
 \indent

For any $p \in  \mathbb{M}(0)$ the non-isotropic group dilations $\delta_{\lambda}:\mathbb{M}(0) \to \mathbb{M}(0)$, $\lambda > 0$, is defined by
\begin{align*}
    \delta_{\lambda}(p_{1},p_{2},p_{3}) = \left(\lambda p_{1}, \lambda p_{2}, \lambda^{2} p_{3}\right).
\end{align*}
Let $\lambda = e^{-t} $ and $\tilde{\gamma} (u,t) : =\delta_{\lambda} (\gamma)$. Then, 
\begin{align*}
    \tilde{\gamma} = \tilde{\gamma}^{1}X_{1} +  \tilde{\gamma}^{2}X_{2} + \tilde{\gamma}^{3}T =\tilde{\gamma}^{1}\partial_x +  \tilde{\gamma}^{2}\partial_y + \tilde{\gamma}^{3}\partial_z =e^{-t}\gamma^{1} \partial_x+ e^{-t}\gamma^{2} \partial_y+ e^{-2t}\gamma^{3}\partial_z 
\end{align*}
and
\begin{align*}
    \tilde{\gamma}^{\prime}(u,t)=\tilde{\gamma}_u^1 \partial_x+\tilde{\gamma}_u^2 \partial_y+\tilde{\gamma}_u^3 \partial_z =\tilde{\gamma}_u^1 X_1+\tilde{\gamma}_u^2 X_2+\left(\tilde{\gamma}_u^3+\tilde{\gamma}^1 \tilde{\gamma}_u^2-\tilde{\gamma}^2 \tilde{\gamma}_u^1\right) T.
\end{align*}
Denote $\tilde{s}$ be the arc-length parameter of $\tilde{\gamma}$, we have 
\begin{align}\label{7.4}
    \tilde{s}(u,t)=|\tilde{\gamma}^{\prime}(u,t)|= \sqrt{(\tilde{\gamma}_{u}^{1})^{2}+(\tilde{\gamma}_{u}^{2})^{2}} = e^{-t} |\gamma^{\prime}(u,t)|,
\end{align}
and 
\begin{align*}
    \dot{\tilde{\gamma}}(u,t) := \frac{\partial \tilde{\gamma}}{ \partial \tilde{s}}, \qquad \qquad |\dot{\tilde{\gamma}}(u,t)| =1.
\end{align*}
Furthermore, 
\begin{align}
    &\dot{\tilde{\gamma}}(u,t) =  \frac{\partial \tilde{\gamma}}{\partial \tilde{s}}=\tilde{\gamma}_{\tilde{s}}^{1}X_{1}\big|_{\tilde{\gamma}} + \tilde{\gamma}_{\tilde{s}}^{2}X_{2}\big|_{\tilde{\gamma}} = \gamma_{s}^{1}X_{1}\big|_{\gamma} +  \gamma_{s}^{2}X_{2}\big|_{\gamma} ,\label{7.5}\\
    & J\dot{\tilde{\gamma}} = -\tilde{\gamma}_{\tilde{s}}^{2}X_{1}\big|_{\tilde{\gamma}} +  \tilde{\gamma}_{\tilde{s}}^{1}X_{2}\big|_{\tilde{\gamma}} =  -\gamma_{s}^{2}X_{1}\big|_{\gamma} +  \gamma_{s}^{1}X_{2}\big|_{\gamma} = J\dot{\gamma}.\label{7.6}
\end{align}

\begin{lemma}\label{lem7.5}
    If $\gamma(u,t)$ is a solution of the flow (\ref{1.2}), then $\tilde{\gamma} (u,t) = \delta_{e^{-t}} (\gamma)$ satisfies the following equation
    \begin{align*}
        \frac{\partial \tilde{\gamma}}{\partial t} = -\left(k^{-1} + g_{\theta}(\tilde{\gamma}, J\dot{\tilde{\gamma}})\right) J\dot{\tilde{\gamma}} - 2\left(\int_{0}^{s} k^{-1} ds  + g_{\theta}(\tilde{\gamma}, T) \right)T .
    \end{align*}
\end{lemma}
\noindent {\bf Proof}:  Since $\tilde{\gamma} (u,t)$ can be written as
\begin{align*}
    \tilde{\gamma} =\tilde{\gamma}^{1}\partial_{x} +\tilde{\gamma}^{2}\partial_{y} +\tilde{\gamma}^{3}\partial_{z}= e^{-t}\gamma^1\partial_{x} + e^{-t}\gamma^2\partial_{y}+e^{-2t}\gamma^3\partial_{z},
\end{align*}
we can compute
\begin{align}\label{7.7}
    \frac{\partial \tilde{\gamma}}{\partial t} &= -\tilde{\gamma} + (e^{-t}\gamma_{t}^{1}, e^{-t}\gamma_{t}^{2}, e^{-2t}\gamma_{t}^{3} - e^{-2t}\gamma^{3}) \notag \\
    & = -\tilde{\gamma} +  e^{-t}\gamma_{t}^{1}X_{1} + e^{-t}\gamma_{t}^{2}X_{2} + e^{-2t}(\gamma_{t}^{3} - \gamma^{3} - \gamma^{1}_{t}\gamma^{2} + \gamma^{2}_{t}\gamma^{2}) T. 
\end{align}
Let $\tilde{k}$ be the curvature of $\tilde{\gamma}$. Then, by (\ref{7.5}) and (\ref{7.6}), we have
\begin{align*}
    \tilde{k} &= g_{\theta}\left(\nabla_{\dot{\tilde{\gamma}}}\dot{\tilde{\gamma}}, J\dot{\tilde{\gamma}}\right) =  g_{\theta}\left(\tilde{\gamma}^{1}_{\tilde{s}\tilde{s}}X_{1} + \tilde{\gamma}^{2}_{\tilde{s}\tilde{s}}X_{2}, -\tilde{\gamma}^{2}_{\tilde{s}}X_{1}+\tilde{\gamma}^{1}_{\tilde{s}}X_{2}\right) =\tilde{\gamma}^{1}_{\tilde{s}} \tilde{\gamma}^{2}_{\tilde{s}\tilde{s}} - \tilde{\gamma}^{2}_{\tilde{s}} \tilde{\gamma}^{1}_{\tilde{s}\tilde{s}}.  
\end{align*}
On the other hand, from (\ref{7.4}), we have 
\begin{align}\label{7.8}
    \frac{\partial }{\partial \tilde{s}} = \frac{e^{t}}{|\gamma^{\prime}(u)|}\frac{\partial }{\partial u} = e^{t} \frac{\partial }{\partial s}, \qquad d\tilde{s} = e^{-t} d s.
\end{align}
Thus,
\begin{align}\label{7.9}
    \tilde{k} = \tilde{\gamma}^{1}_{\tilde{s}} \tilde{\gamma}^{2}_{\tilde{s}\tilde{s}} - \tilde{\gamma}^{2}_{\tilde{s}} \tilde{\gamma}^{1}_{\tilde{s}\tilde{s}} = e^{t}\left(\gamma^{1}_{s}\gamma^{2}_{ss}-\gamma^{2}_{s}\gamma^{1}_{ss}\right) =  e^{t} k.
\end{align}
It is already known that
\begin{align*}
    \frac{\partial \gamma}{\partial t} = -k^{-1}J\dot{\gamma} - 2\int_{0}^{s} k^{-1} ds T = \gamma_{t}^{1}X_{1} +  \gamma_{t}^{2}X_{2} + \left(\gamma_{t}^{3} -  \gamma_{t}^{1}\gamma^{2} +\gamma_{t}^{2}\gamma^{1}\right)T,
\end{align*}
and  
\begin{align*}
    &g_{\theta}\left( \frac{\partial \gamma}{\partial t}, T \right) = \gamma_{t}^{3} -  \gamma_{t}^{1}\gamma^{2} +\gamma_{t}^{2}\gamma^{1} = -2\int_{0}^{s} k^{-1} ds, \\
    &g_{\theta}\left( \frac{\partial \gamma}{\partial t}, J\dot{\gamma} \right) = -\gamma_{t}^{1}\gamma_{s}^{2} + \gamma_{t}^{2}\gamma_{s}^{1} =   -k^{-1},\\
    &g_{\theta}\left( \frac{\partial \gamma}{\partial t},\dot{\gamma} \right) = \gamma_{t}^{1}\gamma_{s}^{1} + \gamma_{t}^{2}\gamma_{s}^{2} = 0.
\end{align*}
Then, by (\ref{7.5})-(\ref{7.9}) and the above equalities, we obtain
\begin{align*}
    &g_{\theta}\left(\frac{\partial \tilde{\gamma}}{\partial t}+\tilde{\gamma} ,  J\dot{\tilde{\gamma}}\right) =  e^{-t}g_{\theta}\left( \frac{\partial \gamma}{\partial t}, J\dot{\gamma} \right)  =  -e^{-t}k^{-1} = - \tilde{k}^{-1} \\
    &g_{\theta}\left(\frac{\partial \tilde{\gamma}}{\partial t}+\tilde{\gamma},  T\right) = -\tilde{\gamma}^{3}+e^{-2t}g_{\theta}\left( \frac{\partial \gamma}{\partial t}, T \right) = -\tilde{\gamma}^{3} -2 \int_{0}^{\tilde{s}} \tilde{k}^{-1} d\tilde{s} \\
    &g_{\theta}\left(\frac{\partial \tilde{\gamma}}{\partial t}+\tilde{\gamma} ,  \dot{\tilde{\gamma}}\right) = e^{-t} g_{\theta}\left( \frac{\partial \gamma}{\partial t},\dot{\gamma} \right) = 0.
\end{align*}
Along $\tilde{\gamma}$, $\left\{\dot{\tilde{\gamma}}, J\dot{\tilde{\gamma}}, T\right\}$ forms an orthonormal basis, thus 
\begin{align}\label{7.10}
    \frac{\partial \tilde{\gamma}}{\partial t} &=  - \tilde{k}^{-1}  J\dot{\tilde{\gamma}} - 2 \left(\int_{0}^{\tilde{s}} \tilde{k}^{-1} d\tilde{s} \right) T -\tilde{\gamma}^{3}T-\tilde{\gamma}  \notag \\
    &=  - \left(\tilde{k}^{-1} +  g_{\theta}(\tilde{\gamma}, J\dot{\tilde{\gamma}})\right)  J\dot{\tilde{\gamma}} - 2\left(\int_{0}^{\tilde{s}}  \tilde{k}^{-1} d\tilde{s} + g_{\theta}(\tilde{\gamma}, T)  \right) T - g_{\theta}(\tilde{\gamma}, \dot{\tilde{\gamma}}) \dot{\tilde{\gamma}} .
\end{align}

Next, it is shown that the geometry of the flow (\ref{7.10}) is determined only by the velocity components in the $J\dot{\tilde{\gamma}}$ and $T$ directions, whereas the tangential component affects only the parameterization of this flow, which is proved in a similar way as in \cite{ZK2001}. Let $\bar{\gamma}(u,t):= \tilde{\gamma}(\varphi(u,t),t) $ and  $\varphi : S^{1} \times [0, \infty)  \to S^{1}  $ satisfies the following equation  
\begin{equation*}
    \begin{cases}
        \frac{\partial \varphi}{\partial t} =  \frac{g_{\theta}(\tilde{\gamma}, \dot{\tilde{\gamma}}) }{|\tilde{\gamma}^{\prime}(u,t)|},\\
        \varphi(u, 0) = u.
    \end{cases}
\end{equation*}
With $\tilde{\gamma}$ already known, the above equation can be regarded as an ODE where $u$ is a parameter. From the smooth dependence on a parameter for solutions of ODE's, we deduce the existence of $\varphi$. Then,  
\begin{align*}
    \frac{\partial \bar{\gamma}}{\partial t} =& \frac{\partial \tilde{\gamma}}{\partial t} (\varphi,t)+ \frac{\partial \tilde{\gamma}}{\partial u}(\varphi,t)\frac{\partial \varphi}{\partial t} 
    =\frac{\partial \tilde{\gamma}}{\partial t} +  \frac{g_{\theta}(\tilde{\gamma}, \dot{\tilde{\gamma}}) }{|\tilde{\gamma}^{\prime}(u,t)|} \tilde{\gamma}^{\prime}(u,t) \\
    =& - \left(\tilde{k}^{-1} +  g_{\theta}(\tilde{\gamma}, J\dot{\tilde{\gamma}})\right)  J\dot{\tilde{\gamma}}- 2\left(\int_{0}^{\tilde{s}}  \tilde{k}^{-1} d\tilde{s} + g_{\theta}(\tilde{\gamma}, T)  \right) T.
\end{align*}
This completes the proof.
\hfill${\square}$

\noindent {\bf Proof of Theorem \ref{thm1.6}}: Combining Proposition \ref{prop3.1} and (\ref{7.4}), we have 
\begin{align*}
    \tilde{\mathcal{L}}_{t} = \int_{\tilde{\gamma}} d\tilde{s}  = \int_{S^{1}} |\tilde{\gamma}^{\prime}(u,t)| du = \int_{S^{1}} e^{-t} | \gamma^{\prime}(u,t)| du = e^{-t} \int_{\gamma} ds = e^{-t} \mathcal{L}_{t} = \mathcal{L}_{0} ,
\end{align*}
that is,
\begin{align}\label{7.11}
    \frac{d}{dt}\tilde{\mathcal{L}}_{t} = 0.
\end{align}
Therefore, the rescaled flow (\ref{1.3}) is a kind of length-preserving flow. 

From Proposition \ref{prop3.2}, $g_{\theta}\left(\dot{\gamma}, T \right) = 0$ holds for all $t > 0$, i.e.,
\begin{align*}
    \gamma_u^3+\gamma^1 \gamma_u^2-\gamma^2 \gamma_u^1=0, \qquad  \forall t >0.
\end{align*}
Since
\begin{align*}
    \tilde{\gamma}_{u}^{1} = e^{-t}\gamma_u^1,\quad \tilde{\gamma}_{u}^{2} = e^{-t}\gamma_u^2,\quad \tilde{\gamma}_{u}^{3} = e^{-2t}\gamma_u^3, 
\end{align*}
we have 
\begin{align}\label{7.12}
    \tilde{\gamma}_{u}^{3}+\tilde{\gamma}^{1}\tilde{\gamma}_{u}^{2}-\tilde{\gamma}^{2} \tilde{\gamma}_{u}^{1}= e^{-2t}\left(\gamma_{u}^{3}+\gamma^{1}\gamma_{u}^{2}-\gamma^{2} \gamma_{u}^{1}\right) = 0
\end{align}
which implies $g_{\theta}\left(\dot{\tilde{\gamma}}, T\right) = 0$, i.e., $\tilde{\gamma}(\cdot, t)$ is a Legendrian curve for $t>0$. Moreover, since the length-preserving flow (\ref{1.3}) is obtained by dilating the flow (\ref{1.2}), the long-time existence of $\tilde{\gamma}(u,t)$ follow directly from Theorem \ref{thm1.1}. 

Denote $\tilde{\gamma}^{*}(u,t) := (e^{-t}\gamma_{1}, e^{-t}\gamma_{2})= e^{-t}\gamma^{*}$ be the projection of $\tilde{\gamma}$ onto $x y$-plane. By Lemma \ref{lem7.1}, $\tilde{\gamma}^{*}(u,t)$ and $\tilde{\gamma}(u,t)$ have the same curvature. According to the theory of the flow (\ref{7.3}) established in \cite{A1998,U1991}, the curvature of limit curve $\tilde{\gamma}^{*}(u,\infty)$ is a nonzero constant. Consequently, the curvature of $\tilde{\gamma}(u,\infty)$ is also a nonzero constant. By Lemma \ref{lem2.6}, it then follows that $\tilde{\gamma}(u,\infty)$ is a Euclidean helix with vertical axis. 
\hfill${\square}$

\noindent {\bf Proof of Proposition \ref{prop1.7}}: A direct calculation yields
\begin{align}\label{7.13}
    \frac{\partial }{\partial t} |\tilde{\gamma}^{\prime}(u)| = |\tilde{\gamma}^{\prime}(u)| \left(1 +  \tilde{k}g_{\theta}(\tilde{\gamma}, J\dot{\tilde{\gamma}})\right).
\end{align}
Thus
\begin{align*}
    \frac{d }{d t}\tilde{\mathcal{L}}_{t}  = \int_{S^{1}}   \frac{\partial }{\partial t} |\tilde{\gamma}^{\prime}(u)| du =  \int_{S^{1}}   \left(1 +  \tilde{k}g_{\theta}(\tilde{\gamma}, J\dot{\tilde{\gamma}})\right) |\tilde{\gamma}^{\prime}(u)|du =  \int_{\tilde{\gamma}_{t}}   \left(1 +  \tilde{k}g_{\theta}(\tilde{\gamma}, J\dot{\tilde{\gamma}})\right) ds .
\end{align*}
By (\ref{7.11}), we obtain
\begin{align*}
    \int_{\tilde{\gamma}_{t}}   \left(1 +  \tilde{k}g_{\theta}(\tilde{\gamma}, J\dot{\tilde{\gamma}})\right) d\tilde{s} = 0, \qquad \forall t \geq 0.
\end{align*}
Since $\tilde{\gamma}_{0} = \delta_{1}(\gamma_{0}) =  \gamma_{0}$ and $\tilde{s}(u,0) = s(u,0) $, we have
\begin{align*}
    \int_{\gamma_{0}}   \left(1 +  k g_{\theta}(\gamma, J\dot{\gamma})\right) ds = 0.
\end{align*}
This completes the proof.
\hfill${\square}$

\noindent {\bf Proof of Proposition \ref{prop1.8}}: From (\ref{7.13}), it follows that
\begin{align*}
    \nabla_{t} \dot{\tilde{\gamma}}  = &\nabla_t\left(\frac{1}{\left|\tilde{\gamma}^{\prime}(u)\right|} \tilde{\gamma}^{\prime}(u)\right) \\
    =&\frac{1}{\left|\tilde{\gamma}^{\prime}(u)\right|} \nabla_u\left(\frac{\partial \tilde{\gamma}}{\partial t}\right)+\frac{1}{\left|\tilde{\gamma}^{\prime}(u)\right|} T_{\nabla}\left(\frac{\partial \tilde{\gamma}}{\partial t}, \frac{\partial \tilde{\gamma}}{\partial u}\right)  -\dot{\tilde{\gamma}}\left(1 +  \tilde{k}g_{\theta}(\tilde{\gamma}, J\dot{\tilde{\gamma}})\right)\\
    =& -\partial_{\tilde{s}}\left(\tilde{k}^{-1} + g_{\theta}(\tilde{\gamma}, J\dot{\tilde{\gamma}}) \right)J\dot{\tilde{\gamma}}  -  \left(\tilde{k}^{-1} + g_{\theta}(\tilde{\gamma}, J\dot{\tilde{\gamma}}) \right)\nabla_{\tilde{s}} J\dot{\tilde{\gamma}} - 2\left(\tilde{k}^{-1}+\nabla_{\tilde{s}} g_{\theta}(\tilde{\gamma} ,T)\right)T \\
    &+2\left( \tilde{k}^{-1} + g_{\theta}(\tilde{\gamma}, J\dot{\tilde{\gamma}})\right)T -\dot{\tilde{\gamma}}\left(1 +  \tilde{k}g_{\theta}(\tilde{\gamma}, J\dot{\tilde{\gamma}})\right) \\
    =& -\partial_{\tilde{s}}\left(\tilde{k}^{-1} + g_{\theta}(\tilde{\gamma}, J\dot{\tilde{\gamma}}) \right)J\dot{\tilde{\gamma}} + 2\left( g_{\theta}(\tilde{\gamma}, J\dot{\tilde{\gamma}}) - \partial_{\tilde{s}}g_{\theta}(\tilde{\gamma} ,T)\right)T
\end{align*}
where the last equality utilizes $\nabla_{\tilde{s}} J\dot{\tilde{\gamma}} = J(\nabla_{\dot{\tilde{\gamma}}}\dot{\tilde{\gamma}})= -\tilde{k}\dot{\tilde{\gamma}}$. Moreover, from (\ref{7.12}), we have
\begin{align*}
    \partial_{\tilde{s}}g_{\theta}(\tilde{\gamma} ,T) = \tilde{\gamma}^{3}_{\tilde{s}} = \tilde{\gamma}^{2} \tilde{\gamma}_{\tilde{s}}^{1}-\tilde{\gamma}^{1}\tilde{\gamma}_{\tilde{s}}^{2} = g_{\theta}(\tilde{\gamma}, J\dot{\tilde{\gamma}}) .
\end{align*}
Thus
\begin{align}\label{7.14}
    \nabla_{t} \dot{\tilde{\gamma}}  = -\partial_{\tilde{s}}\left(\tilde{k}^{-1} + g_{\theta}(\tilde{\gamma}, J\dot{\tilde{\gamma}}) \right)J\dot{\tilde{\gamma}} .
\end{align}
In addition, 
\begin{align}\label{7.15}
    \nabla_{t}\nabla_{\tilde{s}} &=  \frac{1}{\left|\tilde{\gamma}^{\prime}(u)\right|}\nabla_{t}\nabla_{u} -   \left(1 +  \tilde{k}g_{\theta}(\tilde{\gamma}, J\dot{\tilde{\gamma}})\right)\nabla_{\tilde{s}} \notag \\
    & = \nabla_{\tilde{s}}\nabla_{t} -   \left(1 +  \tilde{k}g_{\theta}(\tilde{\gamma}, J\dot{\tilde{\gamma}})\right)\nabla_{\tilde{s}}.
\end{align}
Using (\ref{7.14}) and (\ref{7.15}), we have 
\begin{align*}
    \frac{\partial }{\partial t}\tilde{k}^{2} &= 2g_{\theta}\left(\nabla_{t}\nabla_{\tilde{s}} \dot{\tilde{\gamma}}, \nabla_{\dot{\tilde{\gamma}}}\dot{\tilde{\gamma}}\right) =  2g_{\theta}\left(\nabla_{\tilde{s}}\nabla_{t} \dot{\tilde{\gamma}}, \nabla_{\dot{\tilde{\gamma}}}\dot{\tilde{\gamma}} \right) - 2\tilde{k}^{2}\left(1 +  \tilde{k}g_{\theta}(\tilde{\gamma}, J\dot{\tilde{\gamma}})\right)  \\
    & = -2\tilde{k}\partial^{2}_{\tilde{s}}\left(\tilde{k}^{-1} + g_{\theta}(\tilde{\gamma}, J\dot{\tilde{\gamma}}) \right)  - 2\tilde{k}^{2}\left(1 +  \tilde{k}g_{\theta}(\tilde{\gamma}, J\dot{\tilde{\gamma}})\right).
\end{align*}
Then
\begin{align*}
    \frac{\partial }{\partial t}\tilde{k} = -\partial^{2}_{\tilde{s}}\left(\tilde{k}^{-1} + g_{\theta}(\tilde{\gamma}, J\dot{\tilde{\gamma}}) \right) - \tilde{k}\left(1 +  \tilde{k}g_{\theta}(\tilde{\gamma}, J\dot{\tilde{\gamma}})\right).
\end{align*}
By (\ref{7.13}) and the above equality, we have
\begin{align*}
    \frac{\partial }{\partial t} \int_{\tilde{\gamma}_{t}} \tilde{k} d\tilde{s}  &=  \int_{\tilde{\gamma}_{t}} \frac{\partial \tilde{k}}{\partial t} +\tilde{k} \left(1 +  \tilde{k}g_{\theta}(\dot{\tilde{\gamma}}, J\dot{\tilde{\gamma}})\right) d\tilde{s}  \\
    & =\int_{\tilde{\gamma}_{t}} - \partial^{2}_{\tilde{s}}\left(\tilde{k}^{-1} +  g_{\theta}(\dot{\tilde{\gamma}}, J\dot{\tilde{\gamma}})\right) ds  = 0.
\end{align*}
Therefore 
\begin{align*}
    \int_{\tilde{\gamma}_{t}} \tilde{k} d\tilde{s}  = \int_{\tilde{\gamma}_{\infty}} \tilde{k} d\tilde{s} = \int_{\tilde{\gamma}_{0}} \tilde{k} d\tilde{s}  =  \int_{\gamma_{0}} k ds.
\end{align*}
From Theorem 1.2 of \cite{U1991}, we know that if the initial curve $\tilde{\gamma}^{*}_{0}$ is a smooth, closed, immersed curve with positive curvature, then $\tilde{k}^{*}\big|_{\tilde{\gamma}^{*}_{\infty}}$ is constant and $\tilde{\gamma}^{*}_{\infty}$ is a circle. 

By Lemma \ref{lem7.1}, we have $\tilde{k}\big|_{\tilde{\gamma}_{t}} = \tilde{k}^{*}\big|_{\tilde{\gamma}^{*}_{t}}$ and 
\begin{align*}
    \left|\tilde{\gamma}^{\prime}\right|^2=\left|\left(\tilde{\gamma}^*\right)^{\prime}\right|_{\mathbb{R}^{2}}^2=\left(\tilde{\gamma}_u^1\right)^2+\left(\tilde{\gamma}_u^2\right)^2 
\end{align*}
for all $t\geq 0$, which implies that $\tilde{k}\big|_{\tilde{\gamma}_{\infty}}$ is constant and that $\tilde{\gamma}$ and $\tilde{\gamma}^{*}$ share the same arc-length parameter. Thus 
\begin{align*}
    d\tilde{s} = \left|\tilde{\gamma}^{\prime}\right| d u = \left|\left(\tilde{\gamma}^{*}\right)^{\prime}\right|_{\mathbb{R}^{2}} d u = d\tilde{s}^{*},
\end{align*}
and
\begin{align*}
    \int_{{\gamma}_{0}} k ds = \int_{\tilde{\gamma}_{\infty}} \tilde{k} d\tilde{s} = \int_{S^{1}} \tilde{k}\left|\tilde{\gamma}^{\prime}(u,\infty)\right| du =  \int_{S^{1}} \tilde{k}^{*}\left|(\tilde{\gamma}^{*})^{\prime}(u,\infty)\right| du=\tilde{k}^{*} \tilde{\mathcal{L}}^{*}_{\infty}=2\pi,
\end{align*}
where the last equality utilizes the fact that $\tilde{\gamma}^{*}_{\infty}$ is a circle.
\hfill${\square}$

%%%
%%% Acknowledgments
%%%
%%\section*{\bf Acknowledgments}
%%The authors would like to thank Professors X. Yang,  X. Zhang for their guidance and help! This work is supported by NNSF of China(no.11871275; no.12141104).

%------------------------------------------------------------------------------------%


\begin{thebibliography}{20}
    \bibitem{A1998}{\sc Andrews,\ B.}: \textit{Evolving convex curves}, Calc. Var. Partial Differential Equations {\bf 7}(1998), 315-371.

    \bibitem{B2002}{\sc Blair,\ D.\ E.}: \textit{Riemannian Geometry of Contact and Symplectic Manifolds}, Progress in Mathematics, vol. 203. Birkh$\ddot{\rm a} $user Boston Inc., Boston, MA (2002).

    \bibitem{BN2004}{\sc Bray,\ H.,\  Neves,\ A.}: \textit{Classification of prime 3-manifolds with Yamabe invariant greater than $RP^{3}$}, Ann. of Math. {\bf 159}(2004), 407-424.

    \bibitem{BM2008}{\sc Bray,\ H.,\ Miao,\ P.}: \textit{On the capacity of surfaces in manifolds with nonnegative scalar curvature}, Invent. Math. {\bf 172}(2008), no. 3, 459–475.

    \bibitem{CK2008}{\sc Boyer,\ C.,\ Galicki,\ K.}: \textit{Sasakian geometry}, Oxford Mathematical Monographs, Oxford University Press, Oxford (2008).

    \bibitem{BHW2016}{\sc Brendle,\ S.,\ Hung,\ P.\-K.,\ Wang,\ M.\-T.}: \textit{A Minkowski inequality for hypersurfaces in the anti-de Sitter-Schwarzschild manifold}, Commun. Pure Appl. Math. {\bf 69}(1), (2016), 124-144.

    \bibitem{CDPT2007}{\sc Capogna,\ L., \ Danielli,\ D., Pauls,\ S.,\ Tyson,\ J.}: \textit{An introduction to the Heisenberg group and the sub-Riemannian isoperimetric problem}, Progress in Mathematics 259, Birkhuser Verlag, Basel, 2007. 

    \bibitem{CHL2017}{\sc Chiu,\ H.\ L., \ Huang,\ Y.\ C.,\ Lai,\ S.\ H.}: \textit{An application of the moving frame method to integral geometry in the Heisenberg group}, SIGMA. {\bf 13}(2017), 097.

    \bibitem{CFH2018}{\sc Chiu,\ H.\ L., \ Feng,\ X.\ H., Huang,\ Y.\ C.}: \textit{The differential geometry of curves in the Heisenberg groups}, Differ. Geom. Appl. {\bf 56}(2018), 161-172.

    \bibitem{CH2019}{\sc Chiu,\ H.\ L.,\ Ho,\ P.\ T.}: \textit{Global differential geometry of curves in three-dimensional Heisenberg group and CR sphere}, J. Geom. Anal. {\bf 29}(4) (2019),3438-3469.

    \bibitem{SP2009}{\sc Chanillo,\ S.,\ Yang,\ P.}: \textit{Isoperimetric inequalities $\&$ volume comparison theorems on CR manifolds}, Ann. Sc. Norm. Super. Pisa Cl. Sci. {\bf 8}(2) (2009), 279-307.

    \bibitem{CZ2023}{\sc Cui,\ J.,\ Zhao,\ P.}: \textit{Horizontal inverse mean curvature flow in the Heisenberg group}, arXiv preprint arXiv:2306.15469, 2023.

    \bibitem{CHW2024}{\sc Chang,\ S.\ C., \ Han,\ Y.\ B.,\ Wu,\ C.\ T.}: \textit{Legendrian Mean Curvature Flow in $\eta$-Einstein Sasakian Manifolds}, J. Geom. Anal. {\bf 34}(3), (2024), 89.
        
    \bibitem{DG2006}{\sc Dragomir,\ S.,\ Tomassini,\ G.}: \textit{Differential geometry and analysis on CR manifolds}, Boston, MA, USA: Birkh$\ddot{\rm a} $user Boston, 2006.

    \bibitem{DHW2018}{\sc Drugan,\ G., \ He,\ W.,\ Warren,\ M.\ W.}: \textit{Legendrian curve shortening flow in $\mathbb{R}^{3}$}, Commun. Anal. Geom. {\bf 26}(4), (2018), 759-785.

    \bibitem{FG1996}{\sc Falbel,\ E.,\ Gorodski,\ C.}: \textit{Sub-Riemannian homogeneous spaces in dimensions 3 and 4}, Geom. Dedicata. {\bf 62}(3) (1996), 227-252.

    \bibitem{Gr1989}{\sc Grayson,\ M.\ A.}: \textit{The shape of a figure-eight under the curve shortening flow}, Invent. Math. {\bf 96}(1989), no. 1, 177–180.

    \bibitem{G1990}{\sc Gerhardt,\ C.}: \textit{Flow of nonconvex hypersurfaces into spheres}, J. Differ. Geom. {\bf 32}(1990), no.1, 299-314.

    \bibitem{HI2001}{\sc Huisken,\ G.,\ Ilmanen,\ T.}: \textit{The inverse mean curvature flow and the Riemannian Penrose inequality}, J. Differential Geom. {\bf 59}(2001), no.3, 353-437.

    \bibitem{GR2020}{\sc Gir$\tilde{\rm a} $o,\ F.,\ Rodrigues,\ D.}: \textit{Weighted geometric inequalities for hypersurfaces in sub‐static manifolds}, Bull. Lond. Math. Soc. {\bf 52}(1) (2020), 121-136.

    \bibitem{AC2008}{\sc Hurtado,\ A.,\ Rosales,\ C.}: \textit{Area-stationary surfaces inside the sub-Riemannian three-sphere}, Math. Ann. {\bf 340}(3) (2008), 675-708.

    \bibitem{H2019}{\sc Harvie,\ B.}: \textit{Inverse mean curvature flow over non-star-shaped surfaces}, Math. Res. Lett. {\bf 29}(2022), no.4, 1065-1086.

    \bibitem{HK2024}{\sc Huisken,\ G.,\  Koerber,\ T.}: \textit{Inverse mean curvature flow and Ricci-pinched three-manifolds}, J. Reine Angew. Math. {\bf 2024}(814) (2024), 1-8.

    \bibitem{KWWW2022}{\sc Kwong,\ K.\ K., \ Wei,\ Y.,\  Wheeler,\ V.\ -M.}: \textit{On an inverse curvature flow in two-dimensional space forms}, Math. Ann. {\bf 384}(1) (2022), 1-24.

    \bibitem{KW2023}{\sc Kwong,\ K.\ K.,\  Wei,\ Y.}: \textit{Geometric inequalities involving three quantities in warped product manifolds}, Adv. Math. {\bf 430} (2023), 109213.

    \bibitem{LN2015}{\sc Lee,\ D.\ A.,\  Neves,\ A.}: \textit{The Penrose inequality for asymptotically locally hyperbolic spaces with nonpositive mass}, Comm. Math. Phys. {\bf 339}(2015),  327-352.

    \bibitem{M2007}{\sc Moser,\ R.}:  \textit{The inverse mean curvature flow and p-harmonic functions}, J. Eur. Math. Soc. {\bf 9}(2007), no. 1, 77-83.

    \bibitem{PS2021}{\sc Pan,\ S.,\ Sun,\ J.}: \textit{Curve shortening flow in a 3-dimensional pseudohermitian manifold}, Calc. Var. Partial Differential Equations. {\bf 60}(6), (2021) 212.

    \bibitem{PV2024}{\sc Pisante,\ A.,\ Vecchi,\ E.}: \textit{Global weak solutions for the inverse mean curvature flow in the Heisenberg group}, arXiv preprint arXiv:2406.15123, 2024.

    \bibitem{C2012}{\sc Rosales,\ C.}: \textit{Complete stable CMC surfaces with empty singular set in Sasakian sub-Riemannian 3-manifolds}, Calc. Var. Partial Differential Equations. {\bf 43}(3) (2012), 311-345.
    
    \bibitem{S2003}{\sc Smoczyk,\ K.}: \textit{Closed Legendre geodesics in Sasaki manifolds}, New York J. Math. {\bf 9}(41), (2003), 23-47.

    \bibitem{S2016}{\sc Shi,\ Y.}: \textit{The isoperimetric inequality on asymptotically flat manifolds with nonnegative scalar curvature}, Int. Math. Res. Not. {\bf 2016}(22) (2016), 7038–7050.

    \bibitem{T1969}{\sc Tanno,\ S.}: \textit{Sasakian manifolds with constant $\phi$-holomorphic sectional curvature}, Tohoku Mathematical Journal, Second Series.{\bf 21}(2) (1969), 501-507.

    \bibitem{T1989}{\sc Tanno,\ S.}: \textit{Variational problems on contact Riemannian manifolds}, Trans. Amer. Math. Soc. {\bf 314}(1) (1989), 349-379.

    \bibitem{U1990}{\sc Urbas,\ J.}: \textit{On the expansion of starshaped hypersurfaces by symmetric functions of their principal curvatures}, Math. Z. {\bf 205}(1990), no.1, 355-372.

    \bibitem{U1991}{\sc Urbas,\ J.}: \textit{An expansion of convex hypersurfaces}, J. Differ. Geom. {\bf 31}(1) (1991), 91-125.
    
    \bibitem{ZK2001}{\sc Zhu,\ X.\ P.,\ Chou,\ K.\ S.}: \textit{The curve shortening problem}, Chapman and Hall/CRC, 2001.

    
    
\end{thebibliography}
\end{document}